\documentclass[12pt]{article}
\usepackage{amsmath, amsfonts, amscd, latexsym, amsthm, amssymb}
\newtheorem{theor}{\hspace{1cm}{\sc Theorem}}[section]

\newtheorem{sledst}[theor]{\hspace{1cm}Corollary}
\newtheorem{lemma}[theor]{\hspace{1cm}{\sc Lemma}}
\theoremstyle{definition}
\newtheorem{defin}[theor]{\hspace{1cm}{\sc Definition}}

\newcommand{\codim}{\mathop{\rm codim}\nolimits}

\newcommand{\res}{\mathop{\rm res}\nolimits}

\newcommand{\vol}{\mathop{\rm Vol}\nolimits}
\newcommand{\conv}{\mathop{\rm conv}\nolimits}

\newcommand{\dist}{\mathop{\rm dist}\nolimits}
\newcommand{\MV}{\mathop{\rm MV}\nolimits}
\newcommand{\MF}{\mathop{\rm MF}\nolimits}
\newcommand{\CT}{\mathop{\rm MM}\nolimits}
\newcommand{\CB}{\mathop{\rm CB}\nolimits}

\def\R{\mathbb R}

\def\Z{\mathbb Z}

\def\C{\mathbb C}
\def\CC{({\mathbb C}\setminus 0)}

\begin{document}
\hfill\textit{To Vladimir Igorevich Arnold on the occasion of his
70th birthday.}

\begin{center}
{\Large \textbf{Elimination theory and Newton polytopes.}}
\end{center}

\begin{center}
A. Esterov\footnote{Partially supported by RFBR-JSPS-06-01-91063,
RFBR-07-01-00593 and INTAS-05-7805 grants.}, A.
Khovanskii\footnote{Partially supported by OGP 0156833 grant
(Canada).}
\end{center}

\section{Introduction.} \label{intro}

Let $N\subset\C^n$ be an affine algebraic variety, and let
$\pi:\C^n\to\C^m$ be a projection. The goal of elimination theory
is to describe the defining equations of $\pi(N)$ in terms of the
defining equations of $N$. We shall study defining equations of
projections in the context of Newton polytopes: suppose that the
variety $N\subset\CC^n$ is defined by equations $f_1=\ldots=f_k=0$
with given Newton polytopes and generic coefficients, and the
projection $\pi(N)\subset\CC^m$ is given by one equation $g=0$.
Under this assumption, we shall describe the Newton polytope and
the leading coefficients of the Laurent polynomial $g$ in terms of
the Newton polytope and the leading coefficients of the Laurent
polynomials $f_1,\ldots,f_k$ (by the leading coefficients we mean
the coefficients of monomials from the boundary of the Newton
polytope).

In this section, we define the equation $g$ of a projection of a
complete intersection $f_1=\ldots=f_k=0$ (Definition \ref{defcp})
and describe its Newton polytope in terms of the Newton polytopes
$\Delta_1,\ldots,\Delta_k$ of the equations $f_1,\ldots,f_k$
(Theorem \ref{newtmix}). Section \ref{sectioncomposite} contains
some facts about the geometry of this polytope. In particular,
this polytope is an increasing function of polytopes
$\Delta_1,\ldots,\Delta_k$ (Theorem \ref{thmonot}) and equals the
mixed fiber polytope of $\Delta_1,\ldots,\Delta_k$ up to a shift
and dilatation (Theorem \ref{lcomposmixfib}). The existence and
other basic properties of mixed fiber polytopes (Definition
\ref{defmixfib}) are proved in Section \ref{sectionmixfib}.
Sections \ref{sectiontruncat} and \ref{sectiondeveloped} are
concerned with computation of leading coefficients of $g$. For
example, one can use Theorems \ref{coefofproj} and \ref{kh2gen} to
compute explicitly the coefficients of monomials, which
correspond to the vertices of the Newton polytope of $g$,
provided that the polytopes $\Delta_1,\ldots,\Delta_k$ satisfy
some condition of general position (Definition \ref{defdevel}).
In Section \ref{sectionother}, we present some other versions of
elimination theory in the context of Newton polytopes, such as
elimination theory for rational and analytic functions.

\noindent \textbf{Elimination theory and Newton polytopes.} Many
important problems, related to Newton polytopes and tropical
geometry, turn out to be special cases of this version of
elimination theory. Let us give some examples of such problems.

\noindent \textsc{1) To compute the number of common roots of
polynomial equations with given Newton polytopes and generic
coefficients.} The answer is given by Kouchnirenko-Bernstein's
formula (see \cite{bernst} or Theorem \ref{thbernst} below).

\noindent \textsc{2) To compute the product of common roots of
polynomial equations with given Newton polytopes and generic
coefficients.} If the Newton polytopes satisfy some conditions of
general position, then the answer is given by Khovanskii's product
formula (see \cite{Kh2} or Theorems \ref{khprod} and \ref{kh2gen}
below).

\noindent \textsc{3) To compute the sum of values of a polynomial
over the common roots of polynomial equations with given Newton
polytopes and generic coefficients.} If the Newton polytopes
satisfy some conditions of general position, then the answer is
given by Gelfond-Khovanskii's formula (see \cite{Kh1} or Theorem
\ref{khsum} below).

\noindent \textsc{4) Implicitization theory:} to compute the
Newton polytope and the leading coefficients of the defining
equation of a hypersurface, parameterized by a polynomial mapping
$\CC^n\to\CC^{n+1}$ with given Newton polytopes and generic
coefficients of the components. The Newton polytope was described
by Sturmfels, Tevelev and Yu (see \cite{sturmf2}).

\noindent \textsc{5) To describe the Newton polytope and the
leading coefficients of a multidimensional resultant.} The Newton
polytope and the absolute values of leading coefficients were
computed by Sturmfels (see \cite{sturmf}).

\noindent \textsc{6) To prove the existence of mixed fiber
polytopes (Definition \ref{defmixfib}).} Existence of mixed fiber
polytopes was predicted in \cite{mcd} and proved in \cite{mcm}.

\noindent \textbf{Problems 1-3 in the context of elimination
theory.} To put problems 1-3 in the context of elimination
theory, consider a Laurent monomial as a projection
$\pi:\CC^n\to\CC$. Then the defining equation of the projection
of a 0-dimensional complete intersection
$\{f_1=\ldots=f_n=0\}=\{z_1,\ldots,z_N\}$ is a polynomial
$g(t)=\prod_i \bigl(t-\pi(z_i)\bigr)$ in one variable.

\begin{lemma} \label{lex1} 1) The length of the (1-dimensional) Newton polytope of $g$
equals the number of common roots of $f_1,\ldots,f_n$.\newline 2)
The constant term of $g$ (which is a leading coefficient in our
terminology) equals the product of the values of the monomial
$-\pi$ over all common roots of $f_1,\ldots,f_n$.
\newline 3) Let $S_m$ be the polynomial of $m$ variables, such that $S_m(\sum_i
x_i, \sum_i x_i^2,\ldots, \sum_i x_i^m)$ equals the $m$-th
elementary symmetric function of independent variables $x_i$. Then
the coefficient of the monomial $t^m$ in the polynomial $g$ equals
$(-1)^{N-m}S_m(p_1,\ldots,p_m)$, where $p_m$ is the sum of the
values of the monomial $\pi^m$ over all common roots of
$f_1,\ldots,f_n$.
\end{lemma}
All these facts are obvious, and we omit the proof. We generalize
this lemma to projections of complete intersections of an
arbitrary dimension: see Theorem \ref{newtmix} and Section
\ref{sectiontruncat}.

Lemma \ref{lex1} implies that Kouchnirenko-Bernstein's formula
(Theorem \ref{thbernst}), Khovanskii's product formula (Theorems
\ref{khprod} and \ref{kh2gen}) and Gelfand-Khovanskii's formula
(Theorem \ref{khsum}) can be seen as explicit formulas for the
Newton polytope of $g$, the leading coefficients of $g$, and all
coefficients of $g$ respectively, provided that the Newton
polytopes of $f_1,\ldots,f_n$ satisfy some condition of general
position. We generalize these observations to projections of
complete intersections of an arbitrary dimension: see Section
\ref{sectiondeveloped}.

\noindent \textbf{Problems 4-6 in the context of elimination
theory.} One can consider implicitization theory as a special
case of elimination theory. Indeed, consider a mapping
$g=(g_0,\ldots,g_k):\CC^n\to\CC^{k+1}$ and a $k$-dimensional
complete intersection $F=\{f_1=\ldots=f_{n-k}=0\}\subset\CC^n$
with $g_0,\ldots,g_k,f_1,\ldots,f_{n-k}$ being Laurent
polynomials on $\CC^n$. Let $\pi$ be the standard projection
$\CC^n\times\CC^{k+1}\to\CC^{k+1}$, and let $y_0,\ldots,y_k$ be
the standard coordinates on $\CC^{k+1}$. Then the defining
equation of the image $g(F)\subset\CC^{k+1}$ equals the defining
equation of the projection
$\pi(\{g_0-y_0=\ldots=g_k-y_k=f_1=\ldots=f_{n-k}=0\})$.

A multidimensional resultant is the "universal" special case of
elimination theory, which is clear from the following version of
the definition of a resultant. Consider polynomials
$g_i(x_1,\ldots,x_k)=\sum\limits_{b\in B_i} c_{b,i}x^b$,
$i=0,\ldots,k,\, B_i\subset\Z^k$ as polynomials $f_i$ in
variables $c_{b,i}$ and $x_j$ with all coefficients equal to 1.
Let $\pi$ be the projection of the domain of the polynomials
$(f_0,\ldots,f_k)$ along the domain of the polynomials
$(g_0,\ldots,g_k)$. Then the defining equation of the projection
$\pi(\{f_0=\ldots=f_k=0\})$ is called the
$(B_0,\ldots,B_k)$-\textit{resultant}. Note that this definition
of the multidimensional resultant is somewhat different from the
classical one if we understand the defining equation of a
projection in the sense of Definition \ref{defcp}, since it is not
always square free. We consider the square free version of
Definition \ref{defcp} in Section \ref{sectionother} (see Theorem
\ref{thsqf}).

Elimination theory, implicitization theory and the theory of
multidimensional resultants are equivalent in the sense that they
can be formulated in terms of each other. Thus, the contents of
this paper can be written in terms of resultants or
implcitization theory. When written in these terms, Theorem
\ref{newtmix} turns into the descriptions of Newton polytopes from
\cite{sturmf} and \cite{sturmf2}, while the facts from Sections
\ref{sectiontruncat} and \ref{sectiondeveloped} give some new
information about the leading coefficients. For example, one can
use Theorems \ref{coefofproj} and \ref{kh2gen} to compare the
signs of the leading coefficients of a multidimensional resultant
(see \cite{ekh2}).

Theory of mixed fiber polytopes turns out to be the
Newton-polyhedral counterpart of elimination theory in the
following sense. Define the \textit{composite polytope} of
polytopes $\Delta_1,\ldots\Delta_k$ as the Newton polytope of a
projection of a complete intersection $f_1=\ldots=f_k=0$,
provided that the Newton polytope of $f_i$ is $\Delta_i$ and the
coefficients of $f_1,\ldots,f_k$ are in general position. Then
Theorems \ref{newtmix} and \ref{lcomposmixfib} imply that the
composite polytope satisfies the definition of the mixed fiber
polytope up to a shift and dilatation, which proves the existence
of mixed fiber polytopes. We omit the details and prefer to give
an independent elementary proof of the existence of mixed fiber
polytopes in Section \ref{sectionmixfib} to make our paper
self-contained (the proof from \cite{mcm} is based on the paper
\cite{mcmvoid}, which has not been published by now). Note that
composite polytopes are more convenient then mixed fiber
polytopes in some sense; for example, they are monotonous
(Theorem \ref{thmonot}).

\noindent \textbf{The composite polynomial.} Denote the Zariski
closure of a set $M$ by $\overline M$. For an algebraic mapping
$f:M\to\CC^n$ of an irreducible algebraic variety $M$, denote the
number of points in the preimage $f^{(-1)}(x)$ of a generic point
$x\in f(M)$ by $m(f)$, provided that this number is finite, and
let $m(f)$ be 0 otherwise. Define a \textit{cycle} $N=\sum_i a_i
N_i$ in $\CC^n$ as a formal linear combination of irreducible
algebraic varieties $N_i\subset\CC^n$ of the same dimension with
integer coefficients $a_i$.
\begin{defin} \label{projmult}
Let $\pi:\CC^n\mapsto\CC^{n-k}$ be an epimorphism of complex
tori. For a cycle $N=\sum_i a_i N_i$ in $\CC^n$, the cycle
$\sum_{i} m(\pi|_{N_i})a_i\overline{\pi(N_i)}$ is called
\textit{the projection} $\pi_*N$ of the cycle $N$.
\end{defin}

Let $f_1,\ldots,f_m$ be Laurent polynomials on $\CC^n$, such that
$\codim \{f_1=\ldots=f_m=0\}=m$. Denote the intersection cycle of
the divisors of the polynomials $f_1,\ldots,f_m$ by
$[f_1=\ldots=f_m=0]$.
\begin{defin} \label{defcp}
Let $f_0,\ldots,f_k$ be Laurent polynomials on $\CC^n$, such that
$\codim \{f_0=\ldots=f_k=0\}=k+1$. The Laurent polynomial
$\pi_{f_0,\ldots,f_k}$ on $\CC^{n-k}$, such that
$[\pi_{f_0,\ldots,f_k}=0]=\pi_*[f_0=\ldots=f_k=0]$, is called
\textit{the composite polynomial} of polynomials $f_0,\ldots,f_k$
with respect to the projection $\pi$.
\end{defin}

The composite polynomial $\pi_{f_0,\ldots,f_k}$ is defined up to
a monomial factor. To describe its Newton polytope, we need
Kouchnirenko-Bernstein's formula for the number of roots of a
system of polynomial equations.

\noindent \textbf{Kouchnirenko-Bernstein's formula.} The set of
all convex bodies in $\R^m$ is a semigroup with respect to the
operation of \textit{Minkowskii summation} $A+B=\{a+b\, |\, a\in
A,\, b\in B\}$.
\begin{defin} \textit{The mixed volume} $\MV_{\mu}$, induced by a
volume form $\mu$ on $\R^m$, is the symmetric
Minkowski-multilinear function of $m$ convex bodies in $\R^m$,
such that $\MV_{\mu}(\Delta,\ldots,\Delta)=\int_{\Delta}\mu$ for
every convex body $\Delta\subset\R^m$. The mixed volume, induced
by the standard volume form, is denoted by $\MV$.\end{defin} The
\textit{restriction} $f|_B$ of a Laurent polynomial
$f(x)=\sum_{a\in \Z^n} c_ax^a$ onto a set $B\subset \Z^n$ is the
polynomial $\sum_{a\in B} c_ax^a$. \textit{The Newton polytope}
$\Delta_f$ of a Laurent polynomial $f$ is the convex hull of the
set $A$ such that $f(x)=\sum_{a\in A} c_ax^a$ and $c_a\ne 0$.
\begin{defin} Laurent polynomials $f_0,\ldots,f_k$ on $\CC^n$ are
said to be \textit{Newton-nondegenerate} if, for any collection
of faces $A_0\subset \Delta_{f_0},\ldots,A_k\subset
\Delta_{f_k}$, such that the sum $A_0+\ldots+A_k$ is at most a
$k$-dimensional face of the sum
$\Delta_{f_0}+\ldots+\Delta_{f_k}$, the restrictions
$f_0|_{A_0},\ldots,f_k|_{A_k}$ have no common zeros in $\CC^n$.
\end{defin} Newton-nondegenerate collections of polynomials form a
dense subset in the space of all collections of polynomials with
given Newton polytopes.
\begin{theor}[Kouchnirenko-Bernstein, \cite{bernst}] \label{thbernst} 1) The number of
common roots of Newton-nondegenerate Laurent polynomials
$f_1,\ldots,f_n$ in $\CC^n$, taking multiplicities into account,
is equal to $n!\MV(\Delta_{f_1},\ldots,\Delta_{f_n})$. \newline 2)
Without the assumption of Newton-nondegeneracy, the number of
isolated common roots of $f_1,\ldots,f_n$ in $\CC^n$, taking
multiplicities into account, is not greater than
$n!\MV(\Delta_{f_1},\ldots,\Delta_{f_n})$.
\end{theor}

\noindent \textbf{The Newton polytope of the composite
polynomial.} The Newton polytope of the polynomial
$\pi_{f_0,\ldots,f_k}$ is uniquely determined up to a shift by
the condition $(*)$ below. This condition is a corollary of
Kouchnirenko-Bernstein's formula, and can be seen as its
generalization (see Lemma \ref{lex1}.1).
\begin{theor} \label{newtmix}
1) Let $\pi^{\times}:\Z^{n-k}\hookrightarrow\Z^n$ be the
inclusion of character lattices, defined by the epimorphism
$\pi:\CC^n\mapsto\CC^{n-k}$. Let $A_0,\ldots,A_k\subset\Z^n$ and
$A\subset\Z^{n-k}$ be the Newton polytopes of polynomials
$f_0,\ldots,f_k$ and $\pi_{f_0,\ldots,f_k}$. Then, for any convex
bodies $B_1,\ldots,B_{n-k-1}\subset\Z^{n-k}$,
$$(n-k)!\MV(A,B_1,\ldots,B_{n-k-1})=$$ $$=n!\MV(A_0,\ldots,A_k,\pi^{\times}B_1,\ldots,\pi^{\times}B_{n-k-1}),\eqno(*)$$
provided that the polynomials $f_0,\ldots,f_k$ are
Newton-nondegenerate.
\newline 2) Without the assumption of Newton-nondegeneracy,
$$(n-k)!\MV(A,B_1,\ldots,B_{n-k-1})\leqslant$$ $$\leqslant n!\MV(A_0,\ldots,A_k,\pi^{\times}B_1,\ldots,\pi^{\times}B_{n-k-1}).\eqno(**)$$
\end{theor}

This theorem gives rise to \textit{"elimination theory for convex
bodies"}, which describes the polytope $A$ in terms of
$A_0,\ldots,A_k$, proceeding from the equality $(*)$, and
estimates it, proceeding from the inequality $(**)$. See Section
\ref{sectioncomposite} for details.

\textsc{Proof.} By continuity and linearity of the mixed volume,
it is enough to prove this theorem under the assumption that
$B_1,\ldots,B_{n-k-1}$ are polytopes with integer vertices. Under
this assumption, consider generic Laurent polynomials
$g_1,\ldots,g_{n-k-1}$ on $\CC^{n-k}$ with the Newton polytopes
$B_1,\ldots,B_{n-k-1}$. Since $\pi_{f_0,\ldots,f_k}$ is not
identically zero, the collection $\pi_{f_0,\ldots,f_k},$
$g_1,\ldots,g_{n-k-1}$ is Newton-nondegenerate. If the collection
$f_0,\ldots,f_k$ is Newton-nondegenerate, then the collection
$f_0,\ldots,f_k,g_1\circ\pi,\ldots,g_{n-k-1}\circ\pi$ is also
Newton-nondegenerate.

By Kouchnirenko\,-Bernstein's formula, the number of solutions of
the systems
$f_0=\ldots=f_k=g_1\circ\pi=\ldots=g_{n-k-1}\circ\pi=0$ and
$\pi_{f_0,\ldots,f_k}=g_1=\ldots=g_{n-k-1}=0$ equal
$n!V(A_0,\ldots,A_k,B_1,\ldots,B_{n-k-1})$ and
$(n-k)!V(A,B_1,\ldots,B_{n-k-1})$ respectively. On the other
hand, the solutions of the second system are the projections of
the solutions of the first one. $\Box$

\section{Elimination theory for convex bodies.} \label{sectioncomposite}
Theorem \ref{newtmix} motivates the following definition, which
gives rise to "elimination theory for convex bodies".

A convex body $B$ in an ${(n-k)}$-dimensional subspace
$L\subset\R^n$ is called \textit{a composite body} of convex
bodies $\Delta_0,\ldots,\Delta_k\subset\R^n$, if the mixed volume
$n!\MV(\Delta_0,\ldots,\Delta_k,B_1,\ldots,B_{n-k-1})$ in $\R^n$
equals the mixed volume $(n-k)!\MV(B,B_1,\ldots,B_{n-k-1})$ in
$L$ for every collection of convex bodies
$B_1,\ldots,B_{n-k-1}\subset L$. See Definition \ref{defcompos}
for details.

For every collection of convex bodies $\Delta_0,\ldots,\Delta_k$,
there exists a unique up to a shift composite body (Theorem
\ref{thcompos}). The existence of composite bodies follows from
the fact that the mixed fiber body of bodies
$\Delta_0,\ldots,\Delta_k$ satisfies the definition of a
composite body up to a shift and dilatation (Definition
\ref{defmixfib} and Theorem \ref{lcomposmixfib}).

Thus, the theory of composite bodies is a version of the theory of
mixed fiber polytopes, conjectured in \cite{mcd} and constructed
in \cite{mcm}. Since \cite{mcm} is based on the paper
\cite{mcmvoid}, which have not been published yet, we prefer to
present another approach to mixed fiber polytopes in Section
\ref{sectionmixfib} to make our paper self-contained. At the same
time, we prove some basic facts about composite bodies:\newline -
A composite body of polytopes is a polytope (Theorem
\ref{thmixfib}.2).\newline - A composite body of integer
polytopes (i. e. polytopes such that all their vertices are
integer lattice points) is a shifted integer polytope (Theorem
\ref{thinteger}). \newline - Composite bodies are monotonous
(Theorem \ref{thmonot}).
\newline - The linear span of a composite body depends on the
linear spans of its arguments (Theorem \ref{thspan}).
\newline - Codimension $m$ faces of a composite polytope depend on
codimension $m$ faces of its arguments (Theorem \ref{thfaces}).
\newline - In particular, vertices of the composite polytope of
polytopes $\Delta_0,\ldots,\Delta_k$ can be expressed in terms of
moments of their $k$-dimensional faces (Theorem \ref{thvertices}).

\textbf{Composite bodies.}
\begin{defin} \label{defcompos}
Let $L\subset\R^n$ be a vector subspace of codimension $k$, let
$\mu$ be a volume form on $\R^n/L$, and let
$\Delta_0,\ldots,\Delta_k$ be convex bodies in $\R^n$. A convex
body $B\subset L$ is called a \textit{composite body} of
$\Delta_0,\ldots,\Delta_k$ in $L$ and is denoted by
$\CB_{\mu}(\Delta_0,\ldots,\Delta_k)$ if, for every collection of
convex bodies $B_1,\ldots,B_{n-k-1}\subset L$,
$$n!\MV_{\mu'\wedge\mu}(\Delta_0,\ldots,\Delta_k,B_1,\ldots,B_{n-k-1})=(n-k)!\MV_{\mu'}(B,B_1,\ldots,B_{n-k-1}),$$
where $\mu'$ is a volume form on $L$.
\end{defin}
\begin{theor}\label{thcompos}$
$\newline 1) For any collection of convex bodies
$\Delta_0,\ldots,\Delta_k\subset\R^n$, there exists a composite
body $\CB_{\mu}(\Delta_0,\ldots,\Delta_k)$. \newline 2) A
composite body $\CB_{\mu}(\Delta_0,\ldots,\Delta_k)$ is unique up
to a shift.\end{theor} \textsc{Proof.} Part 1 follows from an
explicit formula for composite bodies, see Theorem
\ref{lcomposmixfib}. Part 2 follows from monotonicity, see
Theorem \ref{thmonot}. $\Box$

The proof of uniqueness implies that, in Definition
\ref{defcompos}, it is enough to consider collections
$B_1,\ldots,B_{n-k-1}$ such that $B_1,\ldots,B_{n-k-1}$ are
simplices. Since composite bodies are unique up to a shift, all
the statements about composite bodies are implied to be valid up
to a shift of a composite body.

\textbf{Monotonicity of a composite body.}
\begin{theor} \label{thmonot}
If $\Delta_i\subset\Delta'_i$ for $i=0,\ldots,k$, then
$\CB_{\mu}(\Delta_0,\ldots,\Delta_k)\subset
\CB_{\mu}(\Delta'_0,\ldots,\Delta'_k)$.
\end{theor}
This is a corollary of monotonicity of mixed volume and the
following fact.
\begin{lemma} Let $\Delta$ and $\Delta'$ be convex bodies in $\R^m$. Suppose that
$$\MV(\Delta,B,\ldots,B)\leqslant\MV(\Delta',B,\ldots,B)$$ for
every simplex $B$. Then, for some shift $a\in\R^n$, the shifted
body $\Delta+a$ is contained in $\Delta'$.
\end{lemma}
\textsc{Proof.} Choose $a$ such that the minimax distance
$$\dist(\Delta+a,\Delta')=\max_{x\in\Delta+a}\min_{y\in\Delta'}|x-y|$$
is minimal. Suppose that $\dist(\Delta+a,\Delta')>0$. Then the set
of all covectors $\gamma\in(\R^m)^*$, such that
$$\max_{x\in\Delta+a}\langle\gamma,x\rangle>\max_{y\in\Delta'}\langle\gamma,y\rangle,$$
is not contained in a half-space. In particular, it contains
covectors $\gamma_0,\ldots,\gamma_m$ such that none of them is a
linear combination of the others with non-negative coefficients.
Denote by $B$ an $m$-dimensional simplex with external normal
covectors $\gamma_0,\ldots,\gamma_m$. Then
$\MV(\Delta,B,\ldots,B)>\MV(\Delta',B,\ldots,B)$ because of the
following formula for mixed volumes. $\Box$
\begin{lemma} Let $\Delta$ be a convex body, let
$B_1,\ldots,B_{m-1}$ be polytopes, and let $\mu$ be a volume form
in $\R^m$. Let $\Gamma\subset(\R^n)^*$ be a set that contains one
external normal covector for each $(m-1)$-dimensional face of the
sum $B_1+\ldots+B_{m-1}$. Then
$$\MV_{\mu}(\Delta,B_1,\ldots,B_{m-1})=\frac{1}{m}\sum_{\gamma\in\Gamma}\max_{x\in\Delta}\langle\gamma,x\rangle
\MV_{\mu/\gamma}(B_1^{\gamma},\ldots,B_{m-1}^{\gamma}),$$ where
$B_i^{\gamma}$ is the maximal face of $B_i$ on which $\gamma$
attains its maximum as a function on $B_i$.
\end{lemma} The mixed volume in the right hand side makes sense, since its arguments
are all parallel to the same $(m-1)$-dimensional subspace
$\ker\gamma$.

\textsc{Proof.} If $\Delta=B_1=\ldots=B_{m-1}$ contains the
origin, then this formula states that the volume of $\Delta$
equals the sum of volumes of the convex hulls $\conv(\{0\}\cup
F)$, where $F$ runs over all $(m-1)$-dimensional faces of
$\Delta$. In general, the formula follows from this special case
by additivity and continuity of the mixed volume. $\Box$

\textbf{Linear span of a composite body.} We need one more fact
about composite bodies, which, in the context of Newton polytopes,
reflects the fact that elimination of variables preserves
homogeneity of equations. Namely, the following theorem expresses
the linear span of a composite body in terms of linear spans of
its arguments.

For a set $\Delta\subset\R^n$, denote the linear span of all
vectors of the form $a-b$, where $a\in\Delta$ and $b\in\Delta$, by
$\langle\Delta\rangle$. For a subspace $L\subset\R^n$, denote the
projection $\R^n\mapsto\R^n/L$ by $p$. Recall that $\mu$ is a
volume form on $\R^n/L$.
\begin{theor}\label{thspan}
1) If $\dim p(\Delta_{i_1}+\ldots+\Delta_{i_q})<q-1$ for some
numbers $0\leqslant i_1<\ldots<i_q\leqslant k$, then
$\CB_{\mu}(\Delta_0,\ldots,\Delta_k)$ consists of one point.
\newline 2) Otherwise, there exists a unique minimal non-empty set
$\{i_1,\ldots,i_q\}\subset\{1,\ldots,n\}$ such that $\dim p
(\Delta_{i_1}+\ldots+\Delta_{i_q})=q-1$. In this case
$$\langle
\CB_{\mu}(\Delta_0,\ldots,\Delta_k)\rangle=\langle\Delta_{i_1}+\ldots+\Delta_{i_q}\rangle\cap
L.$$
\end{theor}
\textsc{Proof.} By definition of a composite body, this theorem
follows from a similar fact about mixed volumes, namely, from D.
Bernstein's criterion for vanishing of the mixed volume (see
below). The uniqueness of a minimal non-empty set
$\{i_1,\ldots,i_q\}\subset\{1,\ldots,n\}$, such that $\dim p
(\Delta_{i_1}+\ldots+\Delta_{i_q})=q-1$, follows from the fact
that the family of all such sets is closed under the operation of
intersection (see \cite{sturmf}, Theorem 1.1 for details). $\Box$

\begin{lemma}[D. Bernstein's criterion, \cite{khovpp}]\label{lbernst} The mixed volume of convex
bodies $B_1,\ldots,B_n$ in $\R^n$ is equal to 0 iff $\dim\langle
B_{i_1}+\ldots+B_{i_q}\rangle<q$ for some numbers $1\leqslant
i_1<\ldots<i_q\leqslant n$.
\end{lemma}

\textbf{Mixed fiber bodies and existence of composite bodies.} The
notion of a composite body turns out to be a version of the
notion of a mixed fiber body. We use this relation to prove the
existence and some basic properties of composite bodies. Recall
the definition of a mixed fiber body.

Let $L\subset\R^n$ be a vector subspace of codimension $k$, let
$\mu$ be a volume form on $\R^n/L$, denote by $p$ the projection
$\R^n\to\R^n/L$.
\begin{defin}[\cite{bs}] For a convex body $\Delta\subset\R^n$, the set of all points of the form
$\int_{p(\Delta)}s\mu\in\R^n$, where $s:p(\Delta)\to\Delta$ is a
continuous section of the projection $p$, is called the
\textit{Minkowski integral} of $\Delta$ and is denoted by $\int
p|_{\Delta}\mu$.
\end{defin}

The following fact explains the relation between composite bodies
and Minkowski integrals.
\begin{lemma} \label{lcomposfib} The convex body
$$(k+1)!{\textstyle\int}p|_{\Delta}\mu$$ is contained in a fiber of the projection $p$ and, up to a shift, satisfies
the definition of the composite body
$\CB_{\mu}(\Delta,\ldots,\Delta)$.
\end{lemma}
\textsc{Proof.} 1) If $\Delta=A+B$, where $B\subset L$ and the
restriction $p|_A$ is injective, then the statement follows from
the additivity of the mixed volume. Indeed, for arbitrary convex
bodies $B_1,\ldots,B_{n-k-1}\subset L$ and a volume form $\mu'$
on $L$,
$$n!\MV_{\mu\wedge\mu'}(A+B,\ldots,A+B,B_1,\ldots,B_{n-k-1})=$$
$$=(k+1)\cdot n!\MV_{\mu\wedge\mu'}(A,\ldots,A,B,B_1,\ldots,B_{n-k-1})=$$ $$=(k+1)\cdot k!\bigl({\textstyle\int}_{p(A)}\mu\bigr)\cdot (n-k)!\MV_{\mu'}(B,B_1,\ldots,B_{n-k-1})=$$
$$=(n-k)!\MV_{\mu'}\Bigl((k+1)!\bigl({\textstyle\int}_{p(A)}\mu\bigr)\cdot
B,B_1,\ldots,B_{n-k-1}\Bigr)=$$
$$=(n-k)!\MV_{\mu'}\Bigl((k+1)!\cdot{\textstyle\int}p|_{\Delta}\mu,B_1,\ldots,B_{n-k-1}\Bigr).$$

2) In general, one can subdivide the projection $p(\Delta)$ into
small pieces, and subdivide $\Delta$ into the inverse images
$\Delta_i$ of these pieces. Representing the mixed volume
$\MV_{\mu\wedge\mu'}(\Delta,\ldots,\Delta,B_1,\ldots,B_{n-k-1})$
as the sum of mixed volumes
$\sum_i\MV_{\mu\wedge\mu'}(\Delta_i,\ldots,\Delta_i,B_1,\ldots,B_{n-k-1})$
for arbitrary convex bodies $B_1,\ldots,B_{n-k-1}$ in $L$, and
approximating each $\Delta_i$ by a sum $A_i+B_i$, such that
$B_i\subset L$ and the restriction $p|_{A_i}$ is injective, one
can reduce the general case to the special case (1). $\Box$

The following theorem provides a way to generalize Lemma
\ref{lcomposfib} to composite bodies of arbitrary collections of
convex bodies.
\begin{theor} \label{thmixfib} Choose a linear projection $u:\R^n\to L$. 1) There exists a unique symmetric multilinear
mapping $\MF_{\mu,u}$ from collections of $k+1$ convex bodies in
$\R^n$ to convex bodies in $L$, such that
$\MF_{\mu,u}(\Delta,\ldots,\Delta)=u\int p|_{\Delta}\mu$ for each
convex body $\Delta\subset\R^n$.\newline 2) This mapping assigns
polytopes to polytopes.
\end{theor} Proof of this theorem is given below.
\begin{defin} \label{defmixfib} The convex body
$\MF_{\mu,u}(\Delta_0,\ldots,\Delta_k)$ is called the
\textit{mixed fiber body} of bodies $\Delta_0,\ldots,\Delta_k$.
\end{defin}
\begin{theor} \label{lcomposmixfib} The convex body
$$(k+1)!\MF_{\mu,u}(\Delta_0,\ldots,\Delta_k)$$ is contained in a fiber of the projection $p$ and, up to a shift, satisfies the
definition of the composite body
$\CB_{\mu}(\Delta_0,\ldots,\Delta_k)$.\end{theor} \textsc{Proof.}
By additivity of mixed fiber bodies and mixed volumes, one can
reduce the statement to the special case
$\Delta_0=\ldots=\Delta_k$ considered in Lemma \ref{lcomposfib}.
$\Box$

\textbf{Virtual bodies.} It is more convenient to prove Theorem
\ref{thmixfib} in the context of virtual bodies instead of convex
bodies, because an explicit formula for mixed fiber bodies (see
Lemma \ref{lexplic}) involves subtraction of convex bodies.

Recall that \textit{the Grothendieck group} $K_G$ of a commutative
semigroup $K$ is the group of formal differences of elements from
$K$. In more details, it is the quotient of the set $K\times K$
by the equivalence relation $(a,b)\sim(c,d)\Leftrightarrow
\exists k : a+d+k=b+c+k$ with operations
$(a-b)+(c-d)=(a+c)-(b+d)$ and $-(a,b)=(b,a)$. For each semigroup
$K$ with cancellation law $a+c=b+c\Rightarrow a=b$, the mapping
$a\to(a+a,a)$ induces the inclusion $K\hookrightarrow K_G$. An
element of the form $(a+a,a)\in K_G$ is said to be
\textit{proper} and is usually identified with $a\in K$. Under
this convention, one can write $(a,b)=a-b$.

\begin{defin} \textit{The group
of virtual bodies} in $\R^n$ is the Grothendieck group of the
semigroup of convex bodies in $\R^n$ with the operation of
Minkowski summation. It contains \textit{the group of virtual
polytopes} in $\R^n$, i. e. the Grothendieck group of the
semigroup of convex polytopes in $\R^n$.
\end{defin} These commutative groups are real vector spaces with the operation of scalar multiplication defined as dilatation.
\begin{defin} For a virtual body
$\Delta$ in $\R^n$, its \textit{support function}
$\Delta(\cdot):(\R^n)^*\to\R$ is defined as
$$\Delta(\gamma)=\max_{x\in\Delta_1}\langle\gamma,x\rangle-\max_{x\in\Delta_2}\langle\gamma,x\rangle,$$
where $\Delta_1$ and $\Delta_2$ are convex bodies such that
$\Delta=\Delta_1-\Delta_2$.
\end{defin}
The following statement describes the group of virtual bodies more
explicitely. A function $f:\R^n\to\R$ is said to be
\textit{positively homogeneous} if $f(tx)=tf(x)$ for each
$t\geqslant 0$. A function $f:\R^n\to\R$ is called \textit{a d.
c. function} if it can be represented as the difference of two
convex functions.
\begin{lemma} 1) The mapping $\Delta\to\Delta(\cdot)$ induces the
isomorphism between the group of virtual bodies in $\R^n$ and the
group of positively homogenious d. c. functions on
$(\R^n)^*$.\newline 2) This isomorphism induces an isomorphism
between the group of virtual polytopes and the group of continuous
piecewise linear positively homogenious functions.
\end{lemma}
\textsc{Proof.} The mapping $\Delta\to\Delta(\cdot)$ is
surjective by definition of a d. c. function. It is injective,
since a convex body is uniquely determined by its support
function. Part 2 follows from the fact that each continuous
piecewise linear function can be represented as a difference of
two convex piecewise linear functions. $\Box$

The operations of taking the mixed volume, the composite body and
the mixed fiber body can be extended to virtual bodies by
linearity. This extension is unique, but its properties are quite
different. For example, The mixed volume of virtual polytopes is
not monotonous (for example, $\MV(-A,A)>\MV(-A,2A)$ for a convex
polygon $A$) and is not non-degenerate in the sense of Lemma
\ref{lbernst} (for example, $\MV(B-C,2B+2C)=0$ for non-parallel
segments $B$ and $C$ in the plane). As a result, virtual
composite bodies do not satisfy Theorems \ref{thmonot} and
\ref{thspan}.

\textbf{Proof of Theorem \ref{thmixfib}.} The uniqueness and Part
2 are corollaries of the following formula for mixed fiber bodies:
\begin{lemma}\label{lexplic} For any convex bodies $\Delta_0,\ldots,\Delta_k\subset\R^n$,
$$\MF_{\mu,u}(\Delta_0,\ldots,\Delta_k)=\frac{1}{(k+1)!}\sum_{0\leqslant i_1<\ldots<i_q\leqslant k}(-1)^{k+1-q}u{\textstyle\int}p|_{(\Delta_{i_1}+\ldots+\Delta_{i_q})}\, \mu.$$\end{lemma}
\textsc{Proof.} Let $m:A\times\ldots\times A\to B$ be a symmetric
multilinear mapping, where $A$ and $B$ are semigroups. Then
$$m(a_1,\ldots,a_k)=\sum\limits_{0\leqslant i_1<\ldots<i_q\leqslant
k}(-1)^{k+1-q}m(a_{i_1}+\ldots+a_{i_q},\ldots,a_{i_1}+\ldots+a_{i_q}).$$
To prove this formula, open the brackets in the right hand side
by linearity of $m$ and cancel similar terms. $\Box$

\begin{lemma} \label{lmixfibpoly}(see Section \ref{sectionmixfib} or \cite{mcm})
Let $u:\R^n\to L$ be a linear projection and let $\mu$ be a
volume form on $\R^n/L$.
\newline 1) There exists a symmetric multilinear mapping $\MF_{\mu,u}$
from collections of $k+1$ virtual polytopes in $\R^n$ to virtual
polytopes in $L$ such that\linebreak
$\MF_{\mu,u}(\Delta,\ldots,\Delta)=u\int p|_{\Delta}\mu$ for each
convex polytope $\Delta\subset\R^n$.
\newline 2) $\MF_{\mu,u}$ maps convex polytopes to convex polytopes.
\end{lemma}
The existence of mixed fiber bodies can be reduced to this
special case as follows. For arbitrary convex bodies
$\Delta_0,\ldots,\Delta_k$ in $\R^n$, define the virtual body
$\MF_{\mu,u}(\Delta_0,\ldots,\Delta_k)$ as in Lemma \ref{lexplic}
above. It follows from the definition that
\newline 1) $\MF_{\mu,u}$ is symmetric, \newline 2)
$\MF_{\mu,u}(\Delta,\ldots,\Delta)=u\int p|_{\Delta}\mu$ for each
convex body $\Delta\subset\R^n$, \newline 3) $\MF_{\mu,u}$ is
continuous in the sense of the norm $|\Delta|=\max_{\gamma\in
B}|\Delta(\gamma)|$, where $B\in(\R^n)^*$ is a compact
neighborhood of the origin, since the Minkowski integral is
continuous in this sense.

Lemma \ref{lmixfibpoly} implies that $\MF_{\mu,u}$ is multilinear
and preserves convexisty under the assumption that the arguments
are polytopes. Namely, for any virtual polytopes
$\Delta_0,\Delta_0',\Delta_1,\ldots,\Delta_k$, \newline 4)
$\MF_{\mu,u}(\Delta_0+\Delta_0',\Delta_1,\ldots,\Delta_k)=\MF_{\mu,u}(\Delta_0,\ldots,\Delta_k)+\MF_{\mu,u}(\Delta_0',\ldots,\Delta_k),$
\newline 5)
$\MF_{\mu,u}(t\cdot
\Delta_0,\ldots,\Delta_k)=t\cdot\MF_{\mu,u}(\Delta_0,\ldots,\Delta_k),$
\newline 6) $\MF_{\mu,u}(\Delta_0,\ldots,\Delta_k)$ is covex if
$\Delta_0,\ldots,\Delta_k$ are convex.

Approximating arbitrary convex bodies with convex polytopes and
using the continuity of $\MF_{\mu,u}$ (property 3), one can extend
properties 4, 5 and 6 to arbitrary convex bodies. $\Box$

\section{Mixed fiber polytopes.} \label{sectionmixfib}
In this section, we prove the existence of mixed fiber polytopes
(Lemma \ref{lmixfibpoly}). Namely, let $L\subset\R^n$ be a vector
subspace of codimension $k$, let $u:\R^n\to L$ be a linear
projection, and let $\mu$ be a volume form on $\R^n/L$. Denote the
projection $\R^n\to\R^n/L$ by $p$. Then
\newline 1) there exists a symmetric multilinear mapping
$\MF_{\mu,u}$ from collections of\linebreak $k+1$ virtual
polytopes in $\R^n$ to virtual polytopes in $L$ such
that\linebreak $\MF_{\mu,u}(\Delta,\ldots,\Delta)=u\int
p|_{\Delta}\mu$ for each convex polytope $\Delta\subset\R^n$.
\newline 2) $\MF_{\mu,u}$ maps convex polytopes to convex polytopes.
\newline \textsc{Proof} follows from the fact that the Minkowski
integral is a polynomial mapping from the space of virtual
polytopes in $\R^n$ to the space of virtual polytopes in $L$.
Every polynomial mapping of vector spaces gives rise to a certain
symmetric multilinear function, which is called the polarization
of the polynomial. In more details, Part 1 follows from Theorems
\ref{thpolar}, \ref{thzardense} and \ref{thminkpolynom}; Part 2
follows from Corollary \ref{cconv}. $\Box$

\textbf{Polarizations of polynomials on Zariski dense sets.} The
existence of mixed fiber polytopes is a corollary of the
following general construction.
\begin{defin}
A set $A$ in a vector space $W$ is said to be \textit{Zariski
dense} if each finite-dimensional subspace $U\subset W$ is
contained in a finite-dimensional subspace $V\subset W$ such that
$A\cap V$ is Zariski closed in $V$ (i. e. $A\cap V$ is not
contained in a proper algebraic subset of $V$).
\end{defin}
\begin{defin}
A map $f:A\to V$ from a subset $A$ of a vector space $W$ to a
vector space $V$ is said to be (homogeneous) polynomial of degree
$k$ if, for each finite-dimensional subspace $U\subset W$ and for
each linear function $l:V\to\R$, the composition $l\circ
f|_U:A\cap U\to\R$ is a restriction of a (homogeneous) polynomial
of degree at most $k$ on $U$.
\end{defin}
\begin{theor} \label{thpolar} 1) A (homogeneous) polynomial map of degree $k$ on a Zariski dense subset of a
vector space $W$ has a unique extension to a (homogeneous)
polynomial map of degree $k$ on $W$. \newline 2) For a
homogeneous polynomial map $f:W\to V$ of degree $k$, there exists
a unique symmetric multilinear function
$Mf:\underbrace{W\oplus\ldots\oplus W}_k\to V$ such that
$Mf(w,\ldots,w)=f(w)$ for every $w\in W$.
\end{theor}
\begin{defin} The function $Mf$ is called \textit{the polarization} of a
polynomial $f$.
\end{defin}
\textsc{Proof.} 1) Let $A$ be a Zariski dense subset in $W$ and
let $f:A\to V$ be a (homogeneous) polynomial map of degree $k$.
For a subspace $U$ such that $A\cap U$ is Zariski dense in $U$,
there exists a unique (homogeneous) polynomial map $f_U:U\to V$ of
degree $k$ such that $f_U=f$ on $U\cap A$. For any two such
finite-dimensional subspaces $U$ and $U'$, the sum $U+U'$ is
contained in a finite-dimensional subspace subspace $U''$ such
that $U''\cap A$ is Zariski dense. Thus $f_{U''}=f_{U'}$ on $U'$
and $f_{U''}=f_{U}$ on $U$. In particular, $f_U=f_{U'}$ on the
intersection $U\cap U'$. This implies that polynomials $f_U$ glue
up into a mapping $\tilde f:W\to V$ such that $\tilde f=f$ on $A$.
\newline 2) For numbers $t_1,\ldots,t_k$ and vectors $w_1,\ldots,w_k\in W$, the expression
$f(t_1w_1+\ldots+t_kw_k)/k!$ is a homogeneous polynomial as a
function of $t_1,\ldots,f_k$. The coefficient of the monomial
$t_1\ldots t_k$ in this polynomial satisfies the definition of
the polarization $Mf$. $\Box$

We apply polarizations in the following context. Let $V(K)$ be
the space of virtual polytopes in a $k$-dimensional vector space
$K$. Let $A(K)\subset V(K)$ be the set of convex polytopes.
\begin{theor}\label{thzardense} $A(K)$ is a Zariski dense subset of $V(K)$.
\end{theor}
\begin{defin} A polytope $\Delta'\in V(K)$ is said to be
compatible with a polytope $\Delta\in V(K)$ if the support
function $\Delta'(\cdot)$ is linear on every domain of linearity
of $\Delta(\cdot)$.
\end{defin}  Let $V(\Delta)\subset V(K)$ be the
space of all virtual polytopes compatible with $\Delta\in V(K)$.
Theorem \ref{thzardense} is a corollary of the following facts.
\newline 1) For every polytope $\Delta\in V(K)$, the
space $V(\Delta)$ is finite dimensional. Indeed, the space of
piecewise-linear functions with the prescribed domains of
linearity is finite-dimensional.
\newline 2) For every convex polytope $\Delta\in A(K)$ the
intersection $V(\Delta)\cap A(K)$ is Zariski dense in $V(K)$.
\newline 3) Every finite dimensional vector subspace $U\subset V(K)$ is
contained in the space $V(\Delta)$ for some convex polytope
$\Delta\in A(K)$. Indeed, if $U$ is generated by differences
$A_i-B_i$ of convex polytopes $A_i$ and $B_i$, then one can
choose $\Delta=\sum_i A_i+B_i$. $\Box$

\textbf{Minkowski integral is a polynomial.} Let $u:\R^n\to L$ be
a linear projection, let $\mu$ be a volume form on the
$k$-dimensional vector space $\R^n/L$, and let $p$ be the
projection $\R^n\to\R^n/L$.
\begin{theor}\label{thminkpolynom} The Minkowski integral $\mathcal{M}(\Delta)=u\int p|_{\Delta}\mu$
is a homogeneous polynomial mapping $A(\R^n)\to A(L)$ of degree
$k+1$.
\end{theor}
\textsc{Proof.} For a convex polytope $\Delta\in\R^n$, define
$A(\Delta)$ as the set of all convex polytopes, compatible with
$\Delta$. For a convex $k$-dimensional polytope $\Delta$, the
restriction of $\mathcal{M}$ to $A(\Delta)$ is a homogeneous
polynomial mapping of degree $k+1$ because of the following two
facts (the first one follows from the definition of the Mnikowski
integral, and the second one is well-known).
\begin{lemma} \label{lminkmoment} The Minkowski integral $\mathcal{M}(\Delta)$ of a convex $k$-dimensional polytope $\Delta$
consists of one point, and this point equals the projection $u$ of
the first moment $\int_{\Delta}xp^*(\mu)$ of $\Delta$, where $x$
runs over $\Delta$ and $p^*(\mu)$ is the volume form $\mu$ on
$\R^n/L$ lifted to $\Delta$. \end{lemma} \begin{lemma}
\label{lmomentpolynom} The first moment is a homogeneous
polynomial of degree $k+1$ on the space $A(\Delta)$, if $\Delta$
is a convex $k$-dimensional polytope.
\end{lemma}

One can reduce Theorem \ref{thminkpolynom} to $k$-dimensional
polytopes as follows. For a covector $\gamma\in(\R^n)^*$ and a
convex polytope $\Delta\subset\R^n$, let $\Delta^{\gamma}$ be the
maximal face, where $\gamma$ attains its maximum as a function on
$\Delta$.
\begin{lemma} \label{lminkfaces} For every covector $\gamma\in L^*$,
$$(u{\textstyle\int}
p|_{\Delta}\mu)^{\gamma}=\sum_{\delta\in(\R^n)^*\atop
\delta|_L=\gamma}u{\textstyle\int} p|_{\Delta^{\delta}}\mu.$$
\end{lemma}
This equality easily follows from the definition of the Minkowski
integral, and we omit the proof. The sum in the right hand side
makes sense, since it contains finitely many non-zero summands.
Note that Lemmas \ref{lminkmoment} and \ref{lminkfaces} are
similar to Propositions 5.1 and 5.2 from \cite{mcm} respectively.
\begin{lemma}\label{lcompat} The mapping $\mathcal{M}$ preserves compatibility of convex polytopes:
$\mathcal{M}\bigl(A(\Delta)\bigr)\subset
A\bigl(\mathcal{M}(\Delta)\bigr)$.
\end{lemma}
\textsc{Proof.} The integral of a continuous family of convex
functions is a linear function iff every function in the family
is linear. Apply this fact to the following description of the
support function of $\mathcal{M}(\Delta)$. $\Box$

For a convex body $\Delta\subset\R^n$ and a point $a\in\R^n/L$,
denote the convex body $u\Bigl(\Delta\cap
p^{(-1)}(a)\Bigr)\subset L$ by $\Delta_a$; roughly speaking, this
is a fiber of $\Delta$ over the point $a$.
\begin{lemma} The support function of the body $u\bigl(\mathcal{M}(\Delta)\bigr)$
equals the integral of the support functions of bodies $\Delta_a$
over $a\in p(\Delta)$.\end{lemma} This equality easily follows
from the definition of the Minkowski integral, and we omit the
proof.

\textsc{Proof of Theorem \ref{thminkpolynom}.} For a face $B$ of
a polytope $\Delta$, let $\tilde B:A(\Delta)\to A(B)$ be the
mapping which maps every $\Delta'\in A(\Delta)$ to its face
$B'\in A(B)$, such that $B+B'$ is a face of $\Delta+\Delta'$.

For an $n$-dimensional convex polytope $\Delta\in A(\R^n)$, denote
vertices of $\mathcal{M}(\Delta)$ by $a_1,\ldots,a_I$, and denote
$k$-dimensional faces of $\Delta$ by $B_1,\ldots,B_J$. By Lemma
\ref{lcompat}, the points $\tilde
a_1\bigl(\mathcal{M}(\Delta')\bigr),\ldots,\tilde
a_I\bigl(\mathcal{M}(\Delta')\bigr)$ are the vertices of the
polytope $\mathcal{M}(\Delta')$ for every convex polytope
$\Delta'\in A(\Delta)$. By Lemma \ref{lminkfaces}, each vertex
$\tilde a_i\bigl(\mathcal{M}(\Delta')\bigr)$ equals a finite sum
of the Minkowski integrals of $k$-dimensional faces $\tilde
B^j(\Delta')$. By Lemmas \ref{lminkmoment} and
\ref{lmomentpolynom}, the Minkowski integral $\mathcal{M}$ is a
homogeneous polynomial of degree $k+1$ on the image of each
linear mapping $\tilde B_j$. $\Box$

\textbf{Faces and convexity of mixed fiber polytopes.} By
Theorems \ref{thzardense} and \ref{thminkpolynom}, there exists a
unique polarization of the Minkowski integral of a polytope in
$\R^n$ with respect to a volume form $\mu$ on $\R^n/L$. It is
denoted by $\MF_{\mu,u}(\Delta_0,\ldots,\Delta_k)$ and is called
the mixed fiber polytope. To prove that it preserves convexity,
we extend Lemma \ref{lminkfaces} to mixed fiber polytopes as
follows.

For a virtual polytope $\Delta$, which equals the difference of
convex polytopes $A$ and $B$ in $\R^n$, and for a covector
$\gamma\in(\R^n)^*$, the \textit{support face} $\Delta^{\gamma}$
is defined as $A^{\gamma}-B^{\gamma}$.
\begin{theor} \label{thfaces}
For virtual polytopes $\Delta_0,\ldots,\Delta_k\subset\R^n$ and a
covector $\gamma\in L^*$, the face
$\bigl(\MF_{\mu,u}(\Delta_0,\ldots,\Delta_k)\bigr)^{\gamma}$
coincides with the Minkowski sum
$$\sum_{\delta\in(\R^n)^*\atop
\delta|_L=\gamma}\MF_{\mu,u}(\Delta_0^{\delta},\ldots,\Delta_k^{\delta}).$$
\end{theor}
This theorem follows from Lemma \ref{lminkfaces} by linearity of
mixed fiber polytopes.

The length of a one-dimensional Minkowski integral of a convex
polytope $\Delta$ is by definition equal to the volume of
$\Delta$. This fact extends by linearity as follows.
\begin{lemma} Suppose that convex polytopes $\Delta_0,\ldots,\Delta_k$ are all
parallel to a $(k+1)$-dimensional subspace $K\subset\R^n$. Choose
a coordinate $t$ on the line $K\cap L$. Then a mixed fiber body
$\MF_{\mu,u}(\Delta_0,\ldots,\Delta_k)$ is a segment, parallel to
the line $K\cap L$, and its length (in the sense of the
coordinate $t$) equals $\MV_{dt\wedge
p^*\mu}(\Delta_0,\ldots,\Delta_k)$.
\end{lemma}
This mixed volume makes sense because its arguments are all
parallel to the same $(k+1)$-dimensional subspace $K$. The volume
form $dt\wedge p^*\mu$ makes sense on $K$, because
$\ker(p|_K)=K\cap L$.

In particular, a one-dimensional mixed fiber polytope of convex
polytopes is convex. Since, by Theorem \ref{thfaces}, every edge
of a mixed fiber polytope is a sum of one-dimensional mixed fiber
polytopes, every edge of a mixed fiber polytope of convex
polytopes is convex. A polytope with all convex edges is convex.
\begin{sledst}\label{cconv} The mixed fiber polytope of convex polytopes is
convex.
\end{sledst}

\textbf{Vertices and integrality of mixed fiber polytopes.} The
proof of Theorem \ref{thminkpolynom} is based on the fact that
vertices of the Minkowski integral of $\Delta$ can be expressed
in terms of the first moments of faces of $\Delta$. We extend
this fact to mixed fiber polytopes in order to prove their
integrality. To formulate this, we need the polarization of the
first moment, which exists by Lemma \ref{lmomentpolynom}. For
virtual polytopes $\Delta_0,\ldots,\Delta_k$ in $\R^n$, the
subspace $\langle\Delta_0,\ldots,\Delta_k\rangle\subset\R^n$ is
defined as the minimal subspace containing convex polytopes
$B_i^j$ such that $\Delta_i=B_i^0-B_i^1$ up to a shift for
$i=0\ldots,k$.
\begin{lemma} \label{ldefmm}
There exists a unique symmetric multilinear function $\CT_{\mu}$
of $k+1$ convex bodies such that \newline 1) The domain of
$\CT_{\mu}$ consists of all collections of virtual
polytopes\linebreak $\Delta_0,\ldots,\Delta_k\subset\R^n$ such
that $\dim\langle\Delta_0,\ldots,\Delta_k\rangle\leqslant k$;
\newline 2) For each $k$-dimensional convex polytope $\Delta\subset\R^n$,
$$\CT_{\mu}(\Delta,\ldots,\Delta)=
\int_{\Delta}xp^*(\mu),$$ where $x$ runs over $\Delta$ and
$p^*(\mu)$ is the volume form $\mu$ on $\R^n/L$ lifted to
$\Delta$.
\end{lemma}
\begin{defin} \label{defmm} The point $\CT_{\mu}(\Delta_0,\ldots,\Delta_k)\in\R^n$ is called the \textit{mixed
moment} of $\Delta_0,\ldots,\Delta_k$.
\end{defin}
By linearity, Lemma \ref{lminkmoment} extends to mixed fiber
polytopes as follows.
\begin{lemma} \label{lmixfibmoment} If $\dim\langle\Delta_0,\ldots,\Delta_k\rangle\leqslant k$, then the
Mixed fiber polytope
$\MF_{\mu,u}(\Delta_0,\ldots,\Delta_k)\in\R^n$ consists of one
point $u\CT_{\mu}(\Delta_0,\ldots,\Delta_k)\in L$.
\end{lemma}
Lemma \ref{lmixfibmoment} and Theorem \ref{thfaces} give the
following expression for vertices of a mixed fiber polytope.
\begin{theor}\label{thvertices} In the notation of Theorem
\ref{thfaces}, \newline 1) If
$\dim\langle\Delta_0^{\delta},\ldots,\Delta_k^{\delta}\rangle\leqslant
k$ for each covector $\delta\in(\R^n)^*$ such that
$\delta|_L=\gamma$, then the face
$\bigl(\MF_{\mu,u}(\Delta_0,\ldots,\Delta_k)\bigr)^{\gamma}$ is a
vertex of $\MF_{\mu,u}(\Delta_0,\ldots,\Delta_k)$, and this vertex
equals $\sum_{\delta|_L=
\gamma}u\CT_{\mu}(\Delta_0^{\delta},\ldots,\Delta_k^{\delta})$.
\newline 2) Almost all covectors $\gamma\in L^*$ satisfy the
condition of Part 1.
\newline 3) The set of all points of the form
$\sum_{\delta|_L=
\gamma}u\CT_{\mu}(\Delta_0^{\delta},\ldots,\Delta_k^{\delta})$,
where $\gamma\in L^*$ satisfies the condition of Part 1,
coincides with the set of all vertices of
$\MF_{\mu,u}(\Delta_0,\ldots,\Delta_k)$.
\end{theor}
In Part 2, "almost all (co)vectors in a space $V$" means "all
covectors from the complement of a finite union of proper vector
subspaces of $V$".

In particular, since the mixed moment of integer polytopes is a
rational number with the denominator $(k+1)!$, the same is true
for mixed fiber polytopes.
\begin{theor}\label{thinteger} If $\Delta_0,\ldots,\Delta_k$ are integer polytopes
(i. e. their vertices are integer lattice points), $L\subset\R^n$
is a $k$-dimensional rational subspace, $u(\Z^n)=L\cap\Z^n$, and
$\mu$ is the integer volume form on $\R^n/L$ (i. e.
$\int_{\R^n/(L+\Z^n)}\mu=1$), then
${(k+1)!}\MF_{\mu,u}(\Delta_0,\ldots,\Delta_k)$ is an integer
polytope.
\end{theor}

\section{Leading coefficients of a composite polynomial in terms of composite polynomials of fewer
variables.}\label{sectiontruncat} We present some technical facts
in this section about how to compute leading coefficients of a
composite polynomial in terms of composite polynomials of fewer
variables. In the next section, we use these facts to compute
leading coefficients of a composite polynomial
$\pi_{f_0,\ldots,f_k}$ explicitly under the assumption that the
Newton polytopes of polynomials $f_0,\ldots,f_k$ satisfy some
condition of general position.

Recall that, for a covector $\gamma\in(\R^n)^*$ and a convex
polytope $A\subset\R^n$, the polytope $A^{\gamma}$ is defined as
the maximal face of $A$, where $\gamma$ attains its maximum as a
function on $A$.
\begin{defin}
For a covector $\gamma\in(\R^n)^*$ and a Laurent polynomial
$f(x)=\sum_{a\in A} c_ax^a$ on $\CC^n$, the polynomial
$\sum_{a\in A^{\gamma}} c_ax^a$  is called the
\textit{truncation} of $f$ in the direction $\gamma$ and is
denoted $f^{\gamma}$. \end{defin} Theorem \ref{reduceproj}
expresses a truncation of a composite polynomial in terms of
composite polynomials of truncations. Theorem \ref{thdehom}
represents a homogeneous composite polynomial as a composite
polynomial of fewer variables. Since truncations of polynomials
are homogeneous, one can use Theorem \ref{thdehom} to simplify the
answer in the formulation of Theorem \ref{reduceproj}. As a
result, one can express a truncation of a composite polynomial in
terms of composite polynomials of fewer variables.
\begin{defin}
\textit{The vertex coefficients} of a polynomial $f$ are the
coefficients of its monomials which correspond to the vertices of
the Newton polytope $\Delta_f$.
\end{defin}
Since a composite polynomial is unique up to a monomial factor,
we are interested in ratios of its vertex coefficients rather
than in individual vertex coefficients. Theorem \ref{coefofproj}
expresses the ratio of two vertex coefficients of a composite
polynomial as the product of values of some monomial over the
roots of some system of polynomial equations. By Lemma
\ref{lex1}.2, this product over roots can be seen as a vertex
coefficient of a corresponding composite polynomial of one
variable.

\noindent\textbf{Truncation and dehomogenization.} The operations
of truncating and taking the composite polynomial commute in the
following sense.
\begin{theor} \label{reduceproj}
Let $\pi:\CC^n\to\CC^{n-k}$ and
$\pi^{\times}:\Z^{n-k}\hookrightarrow\Z^{n}$ be an epimorphism of
complex tori and the corresponding embedding of their character
lattices, and let $f_0,\ldots,f_k$ be Newton-nondegenerate
Laurent polynomials on $\CC^n$. Then, for every
$\gamma\in(\Z^{n-k})^*$, the truncation
$\pi_{f_0,\ldots,f_k}^{\gamma}$ equals the product
$\prod_{\delta} \pi_{f_0^{\delta},\ldots,f_k^{\delta}}$ over all
$\delta\in(\Z^n)^*$ such that
$\delta|_{\pi^{\times}\Z^{n-k}}=\gamma$.
\end{theor}
Since composite polynomials are defined up to a monomial
multiplier, we can assume that whenever
$\pi_{f_0^{\delta},\ldots,f_k^{\delta}}$ is a monomial, it is
equal to 1. Under this assumption, the product
$\prod\limits_{\delta\in \Z^n,\,
\delta|_{\pi^{\times}\Z^{n-k}}=\gamma}
\pi_{f_0^{\delta},\ldots,f_k^{\delta}}$ contains a finite number
of factors different from 1. The proof of this theorem is given
at the end of this section. Theorem \ref{thfaces} is the
geometrical counterpart of this theorem.

A Laurent polynomial $f:\CC^n\to\C$ is said to be
\textit{homogeneous}, if there exist an epimorphism of complex
tori $\CC^n\to\CC^{n'}$ and a Laurent polynomial
$g:\CC^{n'}\to\C$, such that $n'<n$ and $f=g\circ h$ up to a
monomial factor. The polynomial $g$ is called a
\textit{dehomogenization} of $f$. Theorem \ref{thdehom} below
implies that the operations of dehomogenization and taking  the
composite polynomial commute in the following sense: if
polynomials $f_0,\ldots,f_k$ are "homogeneous enough", then their
composite polynomial is also homogeneous, and its
dehomogenization equals the composite polynomial of
dehomogenizations of $f_0,\ldots,f_k$, raised to some power.

Every pair of tori epimorphisms
$\CC^{n-k}\stackrel{\pi}{\leftarrow}\CC^n\stackrel{h}{\rightarrow}\CC^{n'}$
and corresponding character lattice embeddings
$\Z^{n-k}\stackrel{\pi^{\times}}{\hookrightarrow}\Z^n\stackrel{h^{\times}}{\hookleftarrow}\Z^{n'}$
can be included into the commutative squares
$$
\begin{array}{ccccccc}
  \CC^n & \stackrel{h}{\mapsto} & \CC^{n'} & & \Z^n & \stackrel{h^{\times}}{\hookleftarrow} & \Z^{n'} \\
  \phantom{\pi}\downarrow\pi &  & \phantom{\pi'}\downarrow\pi' & \quad \mbox{ and } \quad & \phantom{\pi^{\times}}\uparrow\pi^{\times} &
  & \phantom{\pi'^{\times}}\uparrow\pi'^{\times} \\
  \CC^{n-k} & \stackrel{h'}{\mapsto} & \CC^{n'-k} & & \Z^{n-k} & \stackrel{h'^{\times}}{\hookleftarrow} &
  \Z^{n'-k},
\end{array}
$$
such that the image of $\Z^{n'-k}$ in $\Z^n$ equals the
intersection $\pi^{\times}\Z^{n-k}\cap h^{\times}\Z^{n'}$.
\begin{theor} \label{thdehom} In this notation, if Laurent polynomials $f_0,\ldots,f_k$ on $\CC^n$
are homogeneous in the sense that $f_i=g_i\circ h$ up to a
monomial factor for some Laurent polynomials $g_0,\ldots,g_k$ on
$\CC^{n'}$, then their composite polynomial
$\pi_{f_0,\ldots,f_k}$ is homogeneous in the sense that it equals
$g\circ h'$, where $g=(\pi'_{g_0,\ldots,g_k})^{\bigl|\Z^n\, /\,
(\pi^{\times}\Z^{n-k}+h^{\times}\Z^{n'})\bigr|}$ is a Laurent
polynomial on $\CC^{n'-k}$.
\end{theor}
The proof is given at the end of this section.

\noindent\textbf{Vertex coeficients.}
\begin{defin} \label{prodroots}
\textit{The product over roots}
$R_{A_1,\ldots,A_m}(g_0;g_1,\ldots,g_m)$ is a raitonal function on
the space of collections of Laurent polynomials
$(g_1,\ldots,g_m)$ such that the Newton polytope of $g_i$ is
$A_i\subset\Z^m$. By definition, this function equals the product
of values of a polynomial $g_0$ over the roots of the system
$g_1=\ldots=g_m=0$ for Newton-nondegenerate polynomials
$g_1,\ldots,g_m$.
\end{defin}
Part 2 of Lemma \ref{lex1} is a formula for the vertex
coefficient of a composite polynomial of one variable in terms of
products over roots. The following theorem extends this formula to
composite polynomials of several variables. Let
$\pi:\C^k\times\C^{n-k}\to\C^{n-k}$ be the standard projection,
and let $u_1,\ldots,u_k$ be the standard coordinates on $\C^k$.
Suppose that $f_0,\ldots,f_k$ are polynomials on
$\C^k\times\C^{n-k}$, and their Newton polytopes
$A_0,\ldots,A_k\subset\Z^k\times\Z^{n-k}$ intersect all coordinate
hyperplanes. Denote the Newton polytope of the composite
polynomial $\pi_{f_0,\ldots,f_k}$ by $A\subset\Z^{n-k}$, and
consider covectors $\gamma_1$ and $\gamma_2$ in $(\Z^{n-k})^*$
with positive integer coordinates. Let $\tilde f_i(u,t)$ be a
Laurent polynomial $f_i(u,\, t^{\gamma_2}+t^{-\gamma_1})$ of $k+1$
variables $u_1,\ldots,u_k,t$, and let $\tilde A_i$ be its Newton
polytope.
\begin{theor} \label{coefofproj} If the polynomials
$f_0,\ldots,f_k$ are Newton-nondegenerate, and covectors
$\gamma_1$ and $\gamma_2$ are generic in the sense that, for every
$a\in\Z^k$, the face
$\Bigl((\{a\}\times\Z^{n-k})\cap\sum_iA_i\Bigr)^{\gamma_j}$ is a
vertex, then
\newline 1) the face $A^{\gamma_j}$ of the
polytope $A$ is a vertex (denote it by $B_j$), and the difference
$B_1-B_2$ equals
$$(k+1)!\sum_{\delta\in(\Z^k)^*}\CT_{\mu}(A_0^{\gamma_1+\delta},\ldots,A_k^{\gamma_1+\delta})-\CT_{\mu}(A_0^{\gamma_2+\delta},\ldots,A_k^{\gamma_2+\delta}),$$
where $\mu$ is the unit volume form on $\Z^k$, and the mixed
moment $\CT$ is defined in Lemma \ref{ldefmm};\newline 2) the
ratio of the coefficients of the composite polynomial
$\pi_{f_0,\ldots,f_k}$ at the vertices $B_1$ and $B_2$ equals
$$(-1)^{\gamma_1\cdot B_1 + \gamma_2\cdot B_2}R_{\tilde
A_0,\ldots,\tilde A_k}(t;\tilde f_0,\ldots,\tilde f_k).$$
\end{theor}
After an appropriate monomial change of coordinates and
multiplication polynomials $f_i$ by appropriate monomials, one can
use this theorem to find the ratio of coefficients of the
composite polynomial $\pi_{f_0,\ldots,f_k}$ at two arbitrary
vertices $B_1$ and $B_2$ of its Newton polytope. If
$\pi_{f_0,\ldots,f_k}$ is homogeneous, then a monomial change of
coordinates is not necessary. If the Newton polytopes of the
polynomials $f_0,\ldots,f_k$ satisfy some condition of general
position (see Definition \ref{defdevel}), then one can use
Theorem \ref{kh2gen} to compute $R(t;\tilde f_0,\ldots,\tilde
f_k)$ explicitly.

\textsc{Proof.} Part 1 follows from Theorem \ref{thvertices}. To
prove Part 2, apply the following lemma to the composite
polynomial $\pi_{f_0,\ldots,f_k}$, multiplied by a monomial in
such a way that its Newton polytope belongs to the positive octant
and intersects all coordinate hyperplanes.

\begin{lemma} Suppose that the Newton polytope $A$ of a polynomial $g$ intersects all
coordinate hyperplanes, and $\gamma_1$ and $\gamma_2$ are
covectors with positive integer components. Then the ratio of the
coefficients of $g$ at the vertices $A^{\gamma_1}$ and
$A^{\gamma_2}$ equals $(-1)^{\gamma_1\cdot B_1+\gamma_2\cdot
B_2}$ times the product of roots of the Laurent poynomial in one
variable $g(t^{\gamma_2}+t^{-\gamma_1})$.
\end{lemma}
This lemma is a corollary of the Vieta theorem.

\noindent\textbf{Proof of Theorem \ref{thdehom}.} Extend the
commutative square
$$
\begin{array}{cccccccc}
  \CC^n & \stackrel{h}{\mapsto} & \CC^{n'} & & \CC^n & \stackrel{p}{\mapsto}\;\; T\phantom{\mapsto\;\;}& \stackrel{p_2}{\mapsto} & \CC^{n'} \\
  \phantom{\pi}\downarrow\pi &  & \phantom{\pi'}\downarrow\pi' & \quad \mbox{ to } \quad & &
  \phantom{p_1}\downarrow {p_1} &  & \phantom{\pi'}\downarrow\pi' \\
  \CC^{n-k} & \stackrel{h'}{\mapsto} & \CC^{n'-k} & & & \CC^{n-k} & \stackrel{h'}{\mapsto} & \CC^{n'-k},
\end{array}
$$
where $p_1$ and $p_2$ are the projections of
$\CC^{n-k}\times\CC^{n'}$ to the multipliers, $T$ is the kernel
of the epimorphism $h'\circ p_1-\pi'\circ
p_2:\CC^{n-k}\times\CC^{n'}\to\CC^{n'-k}$, and
$p=(\pi,h):\CC^n\to\CC^{n-k}\times\CC^{n'}$. The corresponding
commutative diagram of embeddings of character lattices implies
that the image of $p^{\times}$ is a sublattice of index
$q=\bigl|\Z^n\, /\,
(\pi^{\times}\Z^{n-k}+h^{\times}\Z^{n'})\bigr|$ in $\Z^n$. Thus, a
fiber of the epimorphism $p$ consists of $q$ points.

For a cycle $N=\sum_i a_iN_i$ in a complex torus $\CC^m$ and an
epimorphism $p:\CC^n\to\CC^m$, denote the cycle $\sum_i a_i p
^{(-1)}(N_i)$ by $p^{(-1)}(N)$. If $m=n$ and a fiber of $p$
consists of $q$ points, then $p_*\circ p^{(-1)}(N)=q\cdot N$.
Thus, $\pi_*\circ h^{(-1)}=(p_1)_*\circ p_*\circ p^{(-1)}\circ
p_2^{(-1)}=q\cdot(p_1)_*\circ p_2^{(-1)}=q\cdot
h^{(-1)}\circ\pi'_*$. To prove the statement of the theorem, apply
both sides of this equality to the cycle $[g_0=\ldots=g_k=0]$.
$\Box$

\noindent \textbf{Truncations of varieties.} The proof of theorem
\ref{reduceproj} is based on the following definition of a
truncation of a variety (just a more geometric reformulation of
the usual one, see \cite{Kazarn}). By varieties we mean formal
sums of irreducible algebraic varieties of the same dimension with
positive coefficients. By the intersection of varieties we mean
the intersection counting multiplicities, which makes sense for
proper intersections only (the intersection of varieties $V_i$ is
said to be \textit{proper}, if its codimension equals the sum of
codimensions of $V_i$). For an algebraic curve $C\subset\CC^n$,
there exists a unique compactification $\tilde
C=C\sqcup\{p_1,\ldots,p_I\}$ which is smooth near all infinite
points $p_i$. A variety $N\subset\CC^n$ is said to be
$\gamma$-\textit{homogeneous} for a linear function $\gamma$ on
the character lattice of the torus $\CC^n$, if $N$ is invariant
under the action of the corresponding one-parameter subgroup
$\{t^{\gamma}\, |\, t\in\CC\}\subset\CC^n$.
\begin{defin}
1) The truncation of an irreducible curve $C\subset\CC^n$ in the
direction $\gamma\in\Z^n$ is a curve $C^{\gamma}=\sum A_i$, where
the summation is over all infinite points $p_i$ of its
compactification $\tilde C$, and a curve $A_i$ is given by a
parameterization $c_it^{\gamma}$, if $C$ is given by a
parameterization $c_it^{\gamma}+\ldots$ near $p_i$.
\newline 2) The truncation of an arbitrary curve $C=\sum m_iC_i$
in the direction $\gamma\in\Z^n$ is a curve $C^{\gamma}=\sum
m_iC_i^{\gamma}$.
\newline 3) The truncation of an $m$-dimensional
variety $M\subset\CC^n$ in the direction $\gamma\in\Z^n$ is an
$m$-dimensional $\gamma$-homogeneous variety $M^{\gamma}$, such
that for any $\gamma$-homogeneous variety $N$ of dimension $\codim
M + 1$\newline a) if $M^{\gamma}\cap N$ is a curve, then $M\cap
N$ is a curve, and \newline b) under this assumption,
$M^{\gamma}\cap N=(M\cap N)^{\gamma}$.
\end{defin}

\begin{lemma} \label{vartrunk} 1) There exists a unique truncation
of a given variety in a given direction. \newline 2) Let
$f_1=\ldots=f_k=0$ be a Newton-nondegenerate complete
intersection. Then its truncation in a direction $\gamma$ is the
complete intersection $f_1^{\gamma}=\ldots=f_k^{\gamma}=0$.
\newline 3) There is a finite number of different truncations of a given variety.\end{lemma}

\textsc{Proof.} Uniqueness follows from the definition. Existence
is a corollary of the following explicit construction for the
truncation of $M\subset\CC^n$ in the direction $\gamma\in \Z^n$.
Without loss of generality we can assume that
$\gamma=(k,0,\ldots,0)$ and define $M^{\gamma}$ as
$p_1^{-1}(\overline{M}\cap\{x_1=0\})$, wher $x_1,\ldots,x_n$ are
the standard coordinates in $\C^n$, $p_1:\CC^n\to\{x_1=0\}$ is
the standard projection, and
$\overline{M}\subset\C\times\CC^{n-1}$ is the closure of the
variety $M\subset\CC^n\subset\C\times\CC^{n-1}$ counting
multiplicities.

Part 2 also follows from this construction. Indeed, the variety
\linebreak
$p_1^{-1}(\overline{\{f_1^{\gamma}=\ldots=f_k^{\gamma}=0\}}\cap\{x_1=0\})$
is given by the ideal $I$, generated by $\gamma$-truncations of
all elements of the ideal $\langle f_1,\ldots,f_k\rangle$. The
ideal $I$ equals $\langle
f_1^{\gamma},\ldots,f_k^{\gamma}\rangle$, because, for any
relation $\sum g_i f_i^{\gamma}=0$, polynomials $g_i$ are
contained in the ideal $\langle
f_1^{\gamma},\ldots,f_k^{\gamma}\rangle$ (the last fact is
equivalent to vanishing of the first homology group of the Koszul
complex for a regular sequence
$f_1^{\gamma},\ldots,f_k^{\gamma}$).

In general, Part 3 follows from the existance of the $c$-fan, or
the Grobner fan of the ideal of a variety $M$ (see
\cite{Kazarn}). If $M$ is a Newton-nondegenerate complete
intersection, which is the only important case for the proof of
theorem \ref{reduceproj}, then Part 3 follows from Part 2.
Indeed, $\gamma_1$ and $\gamma_2$-truncations of a
Newton-nondegenerate complete intersection $f_1=\ldots=f_k=0$
coincide, if $A^{\gamma_1}=A^{\gamma_2}$, where $A$ is the sum of
the Newton polytopes $\Delta_{f_i}$. $\Box$

\noindent \textbf{Proof of Theorem \ref{reduceproj}.} Theorem
\ref{reduceproj} is a special case of the following fact:
\begin{theor}
Let $\pi:\CC^n\to\CC^{n-k}$ and
$\pi^{\times}:\Z^{n-k}\hookrightarrow\Z^{n}$ be an epimorphism of
complex tori and the corresponding embedding of their character
lattices, and let $M\subset\CC^n$ be a variety. Then the
truncation $\bigl(\pi_*(M)\bigr)^{\gamma}$ equals the sum
$\sum_{\delta} \pi_*(M^{\delta})$ over all $\delta\in(\Z^n)^*$
such that $\delta|_{\pi^{\times}\Z^{n-k}}=\gamma$. (in
particular, there is a finite number of non-empty summands).
\end{theor}

If $M$ is 1-dimensional, then this theorem follows from the
definition of the truncation of a curve. If the dimension is
arbitrary, then the number of non-empty summands is finite by
Lemma \ref{vartrunk}.3, since $M^{\delta_1}=M^{\delta_2},\,
\delta_1\ne\delta_2,\,
\delta_1|_{\pi^{\times}\Z^{n-k}}=\delta_2|_{\pi^{\times}\Z^{n-k}}$
implies $\pi_*M^{\delta_1}=\pi_*M^{\delta_2}=\varnothing$. If a
$(\codim M-k+1)$-dimensional $\gamma$-homogeneous variety
$N\subset\CC^{n-k}$ intersects all summands $\pi_*(M^{\delta})$
properly, then
$$N\cap\sum\limits_{\delta\in (\Z^n)^*,\atop
\delta|_{\pi^{\times}\Z^{n-k}}=\gamma}
\pi(M^{\delta})=\bigl(N\cap\pi(M)\bigr)^{\gamma}\quad
\stackrel{(1)}{\Leftrightarrow}$$

$$\pi_*\Bigl(\sum\limits_{\delta\in (\Z^n)^*,\atop
\delta|_{\pi^{\times}\Z^{n-k}}=\gamma} \pi^{(-1)}(N)\cap
M^{\delta}\Bigr)=\Bigl(\pi_*\bigl(\pi^{(-1)}(N)\cap
M\bigr)\Bigr)^{\gamma}\quad \stackrel{(2)}{\Leftrightarrow}$$

$$\pi_*\Bigl(\sum\limits_{\delta\in (\Z^n)^*,\atop
\delta|_{\pi^{\times}\Z^{n-k}}=\gamma} \pi^{(-1)}(N)\cap
M^{\delta}\Bigr)=\pi_*\Bigl(\sum\limits_{\delta\in (\Z^n)^*,\atop
\delta|_{\pi^{\times}\Z^{n-k}}=\gamma} \bigl(\pi^{(-1)}(N)\cap
M\bigr)^{\delta}\Bigr).$$

Here the last equation follows from the definition of a
truncation of the variety $M$. Equivalence (2) is the statement of
the theorem for a curve $\pi^{(-1)}(N)\cap M$. Equivalence (1) is
a corollary of the following fact:
$\pi_*(A\cap\pi^{(-1)}B)=(\pi_*A) \cap B$ for any varieties
$A\subset\CC^{n}$ and $B\subset\CC^{n-k}$.

\section{Leading coefficients of a composite polynomial: explicit answers for generic Newton
polytopes.}\label{sectiondeveloped}
\begin{defin}
\textit{The edge coefficients} of a polynomial $f$ are the
coefficients of its monomials which correspond to the integer
lattice points on the edges of the Newton polytope $\Delta_f$.
\end{defin}
We can compute explicitly the Newton polytope and the vertex and
edge coefficients of a composite polynomial
$\pi_{f_0,\ldots,f_k}$, provided that the Newton polytopes of the
polynomials $f_0,\ldots,f_k$ satisfy the following condition of
general position.
\begin{defin} \label{defdevel} Polytopes $A_0,\ldots,A_k$ in $\R^n$ are said to be
\textit{developed}, if the following condition is satisfied:
\begin{quotation}
Faces $B_0,\ldots,B_k$ of polytopes $A_0,\ldots,A_k$ sum up to a
$k$-dimensional face of the Minkowski sum $A_0+\ldots+A_k
\Rightarrow B_i$ is a vertex of $A_i$ for some $i$.\end{quotation}
\end{defin}

\noindent\textbf{Elimination theory for polynomials with
developed Newton polytopes.} If the Newton polytopes of
polynomials $f_0,\ldots,f_k$ are developed, then the explicit
computation of the Newton polytope and the vertex and edge
coefficients of the composite polynomial $\pi_{f_0,\ldots,f_k}$
is based on the following facts:

- polynomials $f_0,\ldots,f_k$ are Newton-nondegenerate, and the
assumption of Newton nondegeneracy in Theorems \ref{newtmix}.1,
\ref{reduceproj}, \ref{thdehom} and \ref{coefofproj} is redundant.

- Theorems \ref{reduceproj}, \ref{thdehom} and \ref{coefofproj}
express the vertex and edge coefficients of a composite
polynomial of several variables in terms of composite polynomials
of one variable.

- Passing to the right hand side in the formulation of Theorems
\ref{reduceproj}, \ref{thdehom} and \ref{coefofproj} preserves the
property of Newton polytopes to be developed (see Lemmas
\ref{devpres} and \ref{devpres2} below).

- If $\pi_{f_0,\ldots,f_k}$ is a composite polynomial of one
variable, then Lemma \ref{lex1} implies that Khovanskii's product
formula (Theorems \ref{khprod} and \ref{kh2gen}) and
Gelfand-Khovanskii formula (Theorem \ref{khsum}) can be seen as
explicit formulas for the vertex coefficient and the edge
coefficients of $\pi_{f_0,\ldots,f_k}$ respectively.

\begin{lemma} \label{devpres} 1) In the notation of Theorem \ref{reduceproj}, if the Newton polytopes of
polynomials $f_0,\ldots,f_k$ are developed, then the Newton
polytopes of the polynomials $f_0^{\delta},\ldots,f_k^{\delta}$
are also developed for every covector $\delta$. \newline 2) In the
notation of Theorem \ref{thdehom}, if the Newton polytopes of
polynomials $f_0,\ldots,f_k$ are developed, then the Newton
polytopes of the polynomials $g_0,\ldots,g_k$ are
developed.\end{lemma} These facts follow from definitions, and we
omit the proof.

However, in the notation of Theorem \ref{coefofproj}, the Newton
polytopes of polynomials $\tilde f_0,\ldots,\tilde f_k$ are not
usually developed (regardless of the Newton polytopes of
polynomials $f_0,\ldots,f_k$), and we have to consider the
following (weaker) condition.
\begin{defin} Polytopes $A_1,\ldots,A_n$ in $\R^n$ are said to be
\textit{developed with respect to a point} $b\in\R^n$, if the
following condition is satisfied:
\begin{quotation}
Faces $B_1,\ldots,B_n$ of polytopes $A_1,\ldots,A_n$ sum up to a
face of the Minkowski sum $A_1+\ldots+A_n \Rightarrow B_i$ is a
vertex of $A_i$ for some $i$, unless $B_1+\ldots+B_n$ contains a
segment parallel to the vector $b\in\R^n$.\end{quotation}
\end{defin}
\begin{lemma} \label{devpres2} In the notation of Theorem \ref{coefofproj}, if the Newton polytopes of
polynomials $f_0,\ldots,f_k$ are developed, then the Newton
polytopes of the polynomials $\tilde f_0,\ldots,\tilde f_k$ are
developed with respect to the degree of the monomial $t$.
\end{lemma}

\noindent\textbf{Gelfond-Khovanskii's formula and Khovanskii's
product formula.} \begin{defin} For a collection of polytopes
$A_1,\ldots,A_n$ in $\R^n$, let $\phi_i$ be a non-negative
real-valued function on the boundary $\partial(A_1+\ldots+A_n)$,
such that its zero set is the union of all faces of the form
$B_1+\ldots+B_n$, where $B_1,\ldots,B_n$ are faces of
$A_1,\ldots,A_n$ respectively, and $B_i$ is a vertex. \textit{The
combinatorial coefficient} $C_a$ of a vertex
$a\in(A_1+\ldots+A_n)$ is the local degree of the map
$(\phi_1,\ldots,\phi_n):\partial(A_1+\ldots+A_n)\to\partial\R^n_+$
near $a$, provided that $\phi_1\cdot\ldots\cdot \phi_n=0$ near
$a$.
\end{defin}
In particular, the definition of the combinatorial coefficient
makes sense for all vertices of the sum of developed polytopes.
\begin{defin}
Let $f_1,\ldots,f_n,g$ be Laurent polynomials of variables
$x_1,\ldots,x_n$, and suppose that their Newton polytopes
$A_1,\ldots,A_n$ are developed. \textit{The residue} $\res_{a}
\omega_{f_{\cdot},g}$ of a form $\omega_{f_{\cdot},g}=\frac{g
dx_1\wedge\ldots\wedge dx_n}{f_1\ldots f_n x_1\ldots x_n}$ at a
vertex $a$ of the polytope $\sum A_i$ is defined as the constant
term of the series $g \frac{1}{p(a)} \frac{1}{p/p(a)}$, where $p$
is the product $f_1\cdot\ldots\cdot f_n$, $p(a)$ is its term of
degree $a$ of the polynomial $p$, and $\frac{1}{p/p(a)}$ is the
inverse of the polynomial $p/p(a)$ near the origin.
\end{defin}
\begin{theor}[\cite{Kh1}]\label{khsum}
Let $f_1,\ldots,f_n,h$ be Laurent polynomials on $\CC^n$, and
suppose that their Newton polytopes $A_1,\ldots,A_n$ are
developed. Then the sum of values of $h$ over the roots of the
system $f_1=\ldots=f_n=0$ (multiplicities of the roots being
taken into account) equals $(-1)^n\sum_{a} C_a \res_{a}
\omega_{f_{\cdot},h\det\frac{\partial f_{\cdot}}{\partial
x_{\cdot}}}$, where $a$ runs over all vertices of the polytope
$\sum A_i$.
\end{theor}
Let $\Z_2^{n\times m}$ be the space of $\Z_2$-matrices with $n$
rows and $m$ columns.
\begin{defin}
There exists a unique non-zero function
$\det_2:\Z_2^{n\times(n+1)}\to\Z_2$ which is linear and symmetric
as a function of columns and vanishes at degenerate matrices. It
is called the \textit{2-determinant}.
\end{defin}\begin{defin}
Let $f_1,\ldots,f_n$ be Laurent polynomials on $\CC^n$, and
suppose that their Newton polytopes $A_1,\ldots,A_n$ are
developed. \textit{The Parshin symbol} $[f_1,\ldots,f_n,x^b]_a$ of
the monomial $x^b$ at a vertex $a$ of the polytope $\sum A_i$ is
the product
$$(-1)^{\det_2(a_1,\ldots,a_n,b)}f_1(a_1)^{-\det(b,a_2,\ldots,a_n)}\ldots
f_n(a_n)^{-\det(b,a_1,\ldots,a_{n-1})},$$ where $f_i(a)$ is the
term of degree $a$ of the polynomial $f_i$.
\end{defin}
\begin{theor}[\cite{Kh2}]\label{khprod} Let $f_1,\ldots,f_n$ be Laurent polynomials on $\CC^n$, and
suppose that their Newton polytopes $A_1,\ldots,A_n$ are
developed. Then the product of values of the monomial $x^{A_0}$
over the roots of the system $f_1=\ldots=f_n=0$ (taking
multiplicities of the roots into account) equals $\prod_{a}
[f_1,\ldots,f_n,x^{A_0}]_a^{(-1)^n C_a},$ where $a$ runs over all
vertices of the polytope $\sum A_i$.
\end{theor}
In particular, this product is a monomial as a function of the
vertex coefficients of polynomials $f_1,\ldots,f_n$, and this
theorem can be seen as a multidimensional generalization of the
fact that the constant term of a polynomial in one variable
equals the product of the negatives of its roots.

\noindent\textbf{Khovanskii's product formula for Newton
polytopes, developed with respect to a point.} Lemma
\ref{devpres2} implies, that we have to generalize Theorem
\ref{khprod} to polytopes, developed with respect to a point, to
make it applicable in the context of Theorem \ref{coefofproj}.
For a polytope $A\subset \R^n$ and a concave piecewise-linear
function $v:A\to\R$, denote the polyhedron $\{(a,t) | a\in A,
t\le v(a)\}\subset \R^n\oplus\R^1$ by $N(v)$. Let
$v_1,\ldots,v_n$ be piecewise-linear functions on polytopes
$A_1,\ldots,A_n\subset\R^n$. Denote the union of all bounded
faces of the polytope $\sum_i N(v_i)\subset\R^n\oplus\R^1$ by
$\Gamma$. $\Gamma$ is a topological disc. Let
$\Gamma_j\subset\partial\Gamma$ be the union of all faces that can
be represented as $\sum_i B_i$ where $B_i$ are faces of $N(v_i),
i=1,\ldots,n,$ and $B_j$ is a point. Consider a continuous mapping
$(\phi_1,\ldots,\phi_n):\partial\Gamma\to\R^n_+$, such that the
zero set of a function $\phi_j$ is $\Gamma_j$.
\begin{defin} \label{relvert} Functions $v_1,\ldots,v_n$ are called \textit{developed},
if the image of the mapping $(\phi_1,\ldots,\phi_n)$ is contained
in the boundary of the positive octant $\R^n_+$. A point $a\in
\partial(A_1+\ldots+A_n)$ is called \textit{a vertex of the sum $A_1+\ldots+A_n$
with respect to the functions} $v_j$, if it equals the projection
of some vertex $b\subset\partial\Gamma$ of the sum
$N(v_1)+\ldots+N(v_n)$. In this case $b=b_1+\ldots+b_n$, where
$b_j$ is a vertex of the polyhedron $N(v_j)$, and we denote the
projection of $b_j$ by $a_j$.

\textit{The combinatorial coefficient} $C_a$ of a vertex $a$ is
the local topological degree of the mapping
$(\phi_1,\ldots,\phi_n):\Gamma\to\partial\R^n_+$ at the point
$b$. \textit{The Parshin symbol} $[f_1,\ldots,f_n,x^k]_a$ of the
monomial $x^k$ at this vertex is the product
$$(-1)^{\det_2(a_1,\ldots,a_n,k)}f_1(a_1)^{-\det(k,a_2,\ldots,a_n)}\ldots
f_n(a_n)^{-\det(k,a_1,\ldots,a_{n-1})},$$ where $f_i(a)$ is the
term of degree $a$ of the polynomial $f_i$.
\end{defin}
\begin{theor} \label{kh2gen} If polytopes
$A_1,\ldots,A_n\subset\Z^n$ are developed with respect to $A_0$,
then the function $R_{A_1,\ldots,A_n}(x^{A_0};f_1,\ldots,f_n)$
(see Definition \ref{prodroots}) equals the following monomial in
vertex coefficients of polynomials $f_1,\ldots,f_n$:
$$\prod_{a\;\mbox{\scriptsize is a vertex of }\;A_1+\ldots+A_n\;\mbox{\scriptsize with respect to}\;
v_1,\dots,v_n} [f_1,\ldots,f_n,x^{A_0}]_a^{(-1)^n C_a},$$ where
$v_1,\dots,v_n$ are arbitrary developed functions on the polytopes
$A_1,\ldots,A_n$ such that all vertices of $A_1+\ldots+A_n$ with
respect to $v_1,\dots,v_n$ are integer.
\end{theor}
Theorem \ref{khprod} is a special case of this theorem for
developed polytopes $A_1,\ldots,A_n$ and developed functions
$v_i=0$ on them. The statement of Theorem \ref{kh2gen} is true for
Newton-degenerate polynomials $f_1,\ldots,f_n$, but
$R_{A_1,\ldots,A_n}(x^{A_0};f_1,\ldots,f_n)$ is not equal to the
product of values of the monomial $x^{A_0}$ over the roots of the
system $f_1=\ldots=f_n=0$ in this case.

\textsc{Proof.} The main point in the proof of Theorem
\ref{khprod} (see \cite{Kh2}) is the following fact: if
polynomials $f_1,\ldots,f_n$ depend on a parameter $s\in\CC$, and
their Newton polytopes are developed and do not depend on $s$,
then the product of values of the monomial $x^{A_0}$ over the
roots of the system $f_1=\ldots=f_n=0$ as a function of $s$ is a
monomial, because it is a rational function of $s$ which has no
zeroes and no poles. One can easily verify that the same is true
under the assumption that the Newton polytopes of
$f_1,\ldots,f_n$ are developed with respect to $A_0$, if we
consider the function
$R_{A_1,\ldots,A_n}(x^{A_0};f_1,\ldots,f_n)$ instead of the
product of $x^{A_0}$ over the roots of the system
$f_1=\ldots=f_n=0$. $\Box$

\section{Other versions of elimination theory in the context of Newton
polytopes.}\label{sectionother} In this paper, we discuss common
zeros of Laurent polynomials, the multiplicities of zeros being
taken into account. Of course, one can develop the same theory in
many other contexts. Let us give some examples.

\noindent\textbf{Square free composite polynomials.} This point
of view is usual when discussing Newton polytopes of
multidimensional resultants. For a finite set $A\in\Z^n$, denote
the set of all Laurent polynomials $\sum\limits_{a\in A} c_ax^a$
by $\C[A]$. Consider an epimorphism $\pi:\CC^n\to\CC^{n-k}$ and
finite sets $A_0,\ldots,A_k$ in the character lattice $\Z^n$ of
the complex torus $\CC^n$.

The composite polynomial $\pi_{f_0,\ldots,f_k}$ is not square free
for a collection of polynomials
$(f_0,\ldots,f_k)\in\C[A_0]\oplus\ldots\oplus\C[A_k]$, if the
sets $A_0,\ldots,A_k\subset\Z^n$ are degenerate in some sense.
Let $\pi^0_{f_0,\ldots,f_k}$ be the square free polynomial that
has the same zeros as $\pi_{f_0,\ldots,f_k}$. The theorem stated
below expresses \textit{the square free composite polynomial}
$\pi^0_{f_0,\ldots,f_k}$ in terms of $\pi_{f_0,\ldots,f_k}$.
\begin{defin}
Let $L\subset\Z^n$ be an $(n-k)$-dimensional lattice, and let
$p:\Z^n\to\Z^k$ be the projection along $L$. \textit{The
multiplicity} $d(A_0,\ldots,A_k,L)$ of the collection of finite
sets $A_0,\ldots,A_k\subset\Z^n$ with respect to $L$ is defined
as follows:

1) if $\dim p(A_{i_1}+\ldots+A_{i_q})<q-1$ for some numbers
$0\leqslant i_1<\ldots<i_q\leqslant k$, then
$d(A_0,\ldots,A_k,L)=0$;

2) Otherwise, choose the minimal non-empty set
$\{i_1,\ldots,i_q\}\subset\{0,\ldots,k\}$ such that $\dim
p(A_{i_1}+\ldots+A_{i_q})=q-1$, choose the minimal sublattice
$M\subset\Z^n$ that contains the sum $A_{i_1}+\ldots+A_{i_q}+L$ up
to a shift, and note that $\codim M=k+1-q$. Denote the projection
$\Z^n\to\Z^{k+1-q}$ along $M$ by $r$, and denote the set
$\{0,\ldots,k\}\setminus\{i_1,\ldots,i_q\}$ by
$\{j_1,\ldots,j_{k+1-q}\}$. In this notation,
$$d(A_0,\ldots,A_k,L)=(k+1-q)!\MV(rA_{j_1},\ldots,rA_{j_{k+1-q}})\cdot|\ker
r/M|.$$
\end{defin}
For example, suppose that $L\subset\Z^2$ is the horizontal
coordinate axis. In this case $d(A_1,A_2,L)=0$ iff both $A_1$ and
$A_2$ are contained in horizontal segments. If one of them is
contained in a horizontal segment, then $d(A_1,A_2,L)$ equals the
height of the other one. If neither $A_1$ nor $A_2$ is contained
in a horizontal segment, then $d(A_1,A_2,L)$ equals the GCD of
lengths of vertical segments, connecting points of the set
$A_1+A_2+L$.
\begin{theor} \label{thsqf}
Consider an epimorphism of complex tori $\pi:\CC^n\to\CC^{n-k}$,
the corresponding embedding of their character lattices
$L\subset\Z^n$, and finite sets $A_0,\ldots,A_k\subset\Z^n$.

1) If $d(A_0,\ldots,A_k,L)=0$, then
$\pi^0_{f_0,\ldots,f_k}=\pi_{f_0,\ldots,f_k}=1$ for all
collections of polynomials
$(f_0,\ldots,f_k)\in\C[A_0]\oplus\ldots\oplus\C[A_k]$.

2) Otherwise,
$(\pi^0_{f_0,\ldots,f_k})^{d(A_0,\ldots,A_k,L)}=\pi_{f_0,\ldots,f_k}$
for all collections of polynomials $f_0,\ldots,f_k$ from some
Zariski open subset of the space
$\C[A_0]\oplus\ldots\oplus\C[A_k]$.
\end{theor}
Note that this Zariski open subset neither contains nor is
contained in the set of all Newton-nondegenerate collections of
polynomials. In particular, this theorem implies that the Newton
polytope of $\pi^0_{f_0,\ldots,f_k}$ is $d(A_0,\ldots,A_k,L)$
times smaller, than the Newton polytope of $\pi_{f_0,\ldots,f_k}$
for a generic collection of polynomials
$(f_0,\ldots,f_k)\in\C[A_0]\oplus\ldots\oplus\C[A_k]$.

\textsc{Proof.} Part 1 is a corollary of Theorem \ref{thspan}.1.
Applying Theorem \ref{thspan}.2 and Theorem \ref{thdehom}, one
can reduce Part 2 to the following special case. $\Box$
\begin{defin} Let $L\subset\Z^n$ be an $(n-k)$-dimensional lattice
$L\subset\Z^n$, and let $p:\Z^n\to\Z^k$ be the projection along
$L$. A collection of finite sets $A_0,\ldots,A_k\subset\Z^n$ is
said to be \textit{essential} with respect to $L\subset\Z^n$ if
$\dim p(A_{i_1}+\ldots+A_{i_q})>q-1$ for every collection of
numbers $0\leqslant i_1<\ldots<i_q\leqslant k,\, q\leqslant k$,
and the sum $A_0+\ldots+A_k+L$ is not contained in a shifted
proper sublattice of $\Z^n$.
\end{defin}
Note that $d(A_0,\ldots,A_k,L)=1$ if the collection
$A_0,\ldots,A_k\subset\Z^n$ is essential with respect to $L$, but
the converse is not true. For example, if $L\subset\Z^2$ is the
horizontal coordinate axis, then $A_1$ and $A_2$ form an essential
collection iff neither of them is contained in a horizontal
segment and $d(A_1,A_2,L)=1$.
\begin{lemma} If $A_0,\ldots,A_k$ are essential
with respect to $L$, then
$\pi^0_{f_0,\ldots,f_k}=\pi_{f_0,\ldots,f_k}$ for all collections
of polynomials $f_0,\ldots,f_k$ from some Zariski open subset of
the space $\C[A_0]\oplus\ldots\oplus\C[A_k]$.
\end{lemma}
The proof is the straightforward generalization of a similar
argument for multidimensional resultants (see \cite{sturmf},
Theorem 1.1, where the notion of essential sets was introduced
for $k=n$).

\noindent\textbf{Composite functions of rational functions.}
\begin{defin} \label{virtnewt} \textit{A vertex} of a virtual
polytope $A-B$ is a pair of vertices $(a,b)$ of polytopes $A$ and
$B$, such that $a+b$ is a vertex of the sum $A+B$. \textit{The
Newton polytope of a rational function} $\frac{f}{g}$ is the
difference of the Newton polytopes of $f$ and $g$. \textit{The
vertex coefficient of a rational function} $\frac{f}{g}$ at the
vertex $(a,b)$ of its Newton polytope is the ratio of the vertex
coefficients of polynomials $f$ and $g$ at the vertices $a$ and
$b$ respectively.
\end{defin}
One can readily generalize elimination theory from Laurent
polynomials and convex polytopes to rational functions and virtual
polytopes.

\noindent\textbf{Composite functions of germs of analytic
functions.} A \textit{convex polyhedron} in $\R^n$ is an
intersection of a finite number of half-spaces (which may be
unbounded). Two convex polyhedra in $\R^n$ are said to be
\textit{parallel} if their support functions have the same
domain. For a germ of an analytic function $f:(\C^n,0)\to(\C,0)$
of variables $x_1,\ldots,x_n$, the Newton polyhedron is defined as
the minimal polyhedron parallel to the positive octant in the
lattice of monomials in $x_1,\ldots,x_n$ and containing all
monomials of the Taylor expansion of $f$. One can readily
generalize elimination theory from Laurent polynomials and
bounded polyhedra to germs of analytic functions and polyhedra
parallel to the positive octant. It requires the following
version of Bernstein's theorem.

\begin{defin}[\cite{E1}, \cite{E2}] \label{mixpairs} Let $P_C$ be the set of all pairs of polyhedra $(A,B)$,
such that $A$ and $B$ are both parallel to a cone $C$ and the
difference $A\triangle B$ is bounded. The notions of Minkowski sum
$(A,B)+(C,D)=(A+C,B+D)$ and volume
$\vol\bigl((A,B)\bigr)=\vol(A\setminus B)-\vol(B\setminus A)$ for
such pairs give rise to the mixed volume $V_C:
\underbrace{P_C\times\ldots\times P_C}_{n}\to\R$, which is the
polarization of the volume with respect to Minkowski summation.
\end{defin}

If the cone $C$ consists of one point $*$, then $P_C$ is the set
of pairs of bounded polyhedra,
$\vol\bigl((A,B)\bigr)=\vol(A)-\vol(B)$, and thus
$$V_*((A_1,B_1),\ldots,(A_n,B_n))=\MV(A_1,\ldots,A_n)-\MV(B_1,\ldots,B_n).$$
If $C\ne\{*\}$, then the mixed volumes on the right-hand side are
infinite, but "their difference is well defined".

\begin{lemma}[\cite{E2}] \label{mixelem1}
$V_C((A_1,B_1),\ldots,(A_n,B_n))+$ $$+
V_C((B_1,C_1),\ldots,(B_n,C_n))=V_C((A_1,C_1),\ldots,(A_n,C_n)).$$
\end{lemma}
Let $\mu$ be the unit volume form in $\R^n$, let $S$ be the
positive octant in $(\R^n)^*$, and let $S_0\subset S$ be a set of
covectors that contains a unique multiple of every covector from
$S$.
\begin{lemma}[\cite{E2}] $V_C((A_1,B_1),\ldots,(A_n,B_n))=$ $$=\frac{1}{n} \sum_{\gamma\in S_0}
\sum_{i=1}^n
(\max\gamma(A_i)-\max\gamma(B_i))\MV_{\mu/\gamma}(A^{\gamma}_1,\ldots,A^{\gamma}_{i-1},B^{\gamma}_{i+1},\ldots,B^{\gamma}_n).$$
\end{lemma}
Note that the right hand side of this formula is not symmetric
with respect to permutations of pairs. The sum in the right hand
side makes sense, since it containts finitely many non-zero
summands (which correspond to normal covectors of bounded
$(n-1)$-dimensional faces of the sum $A_1+B_1+\ldots+A_n+B_n$).
The $(n-1)$-dimensional mixed volume in the right hand side makes
sense, since all arguments are contained in the
$(n-1)$-dimensional space $\ker\gamma$.

\begin{defin}
A polyhedron is called an $M$-\textit{far stabilization} of a
polyhedron $\Delta\subset\R^n_+$ parallel to the positive octant
$\R^n_+$, if it can be represented as the convex hull of a union
$\Delta\cup\Gamma$ for some polyhedron $\Gamma\subset\R^n_+$, such
that the distance between $\Gamma$ and the origin is greater than
$M$, and the difference $\R^n_+\setminus\Gamma$ is bounded.
\textit{The mixed volume of (unbounded) polyhedra}
$\Delta_1,\ldots,\Delta_n\subset\R^n_+$ parallel to $\R^n_+$ is
defined as the mixed volume of pairs
$(\R^n_+,\tilde\Delta_1),\ldots,(\R^n_+,\tilde\Delta_n)$, where
$\tilde\Delta_i$ is an $M$-far stabilization of $\Delta_i$,
provided that the mixed volume of these pairs is independent of
the choice of $M$-far stabilizations for some $M$ (we say that
the mixed volume of $\Delta_1,\ldots,\Delta_n$ is
\textit{well-defined} in this case).
\end{defin}

\begin{theor}
1) The mixed volume of polyhedra
$\Delta_1,\ldots,\Delta_n\subset\R^n_+$ parallel to the positive
octant $\R^n_+$ is well-defined if and only if each
$k$-dimensional coordinate plane intersects at least $k$ of these
polyhedra.
\newline 2) If the germs of functions
$f_1,\ldots,f_n:(\C^n,0)\to(\C,0)$ have an isolated common root of
multiplicity $\mu$, then the mixed volume $V$ of their Newton
polyhedra is well-defined, and $\mu\geqslant n!V$.
\newline 3) If the mixed volume $V$ of integer polyhedra
$\Delta_1,\ldots,\Delta_n\subset\R^n_+$ parallel to the positive
octant $\R^n_+$ is well defined, then the germs of analytic
functions $f_1,\ldots,f_n:(\C^n,0)\to(\C,0)$ have an isolated
common root of multiplicity $n!V$, provided that their Newton
polyhedra are equal to $\Delta_1,\ldots,\Delta_n$ and their
leading coefficients are in general position in the following
sense:\begin{quotation} For any collection of bounded faces
$A_1\subset \Delta_1,\ldots,A_n\subset \Delta_n$, such that the
sum $A_1+\ldots+A_n$ is a face of the sum
$\Delta_1+\ldots+\Delta_n$, the Laurent polynomials
$f_1|_{A_1},\ldots,f_n|_{A_n}$ have no common zeros in $\CC^n$.
\end{quotation}
\end{theor}
Part 1 follows from Lemma \ref{mixelem1}, Parts 2 and 3 follow
from Part 1 and a local version of Bernstein's formula (see
\cite{E2}, Theorem 3).

\end{document}